\documentclass[11pt]{article}
\usepackage{xypic}
\usepackage{mathrsfs}
\usepackage{amssymb}
\usepackage{amsfonts}
\usepackage{amsmath}

\def\ben{\begin{enumerate}}
\def\een{\end{enumerate}}

\def\bit{\begin{itemize}}
\def\eit{\end{itemize}}

\def\0{\leqno}

\input xy
\xyoption{all}
\input{tcilatex}
\begin{document}

\begin{center}
{\Large ALGEBRAIC CONSTRUCTIONS }

\ \ \\[0pt]

{\Large IN THE CATEGORY OF VECTOR BUNDLES}

\ \ \\[0pt]


\ \ \\[0pt]
by


\textbf{{CONSTANTIN M. ARCU\c{S} }}

\ \ \\[0pt]

\begin{tabular}{c}
SECONDARY SCHOOL \textquotedblleft CORNELIUS RADU\textquotedblright , \\
RADINESTI VILLAGE, 217196, GORJ COUNTY, ROMANIA \\
e-mail: c\_arcus@yahoo.com, c\_arcus@radinesti.ro%
\end{tabular}
\end{center}



\ \ \\[0pt]





\

\begin{abstract}
The category of generalized Lie algebroids is presented. We obtain an
exterior differential calculus for generalized Lie algebroids. In
particular, we obtain similar results with the classical and modern results
for Lie algebroids. So, a new result of Maurer-Cartan type is presented.
Supposing that any vector subbundle of the pull-back vector bundle of a
generalized Lie algebroid is called \emph{interior differential system (IDS)}
for that generalized Lie algebroid, a theorem of Cartan type is obtained.
Extending the classical notion of \emph{exterior differential system (EDS)}
to generalized Lie algebroids, a theorem of Cartan type is obtained. Using
the theory of linear connections of Ehresmann type presented in the paper $%
\left[ 1\right]$, the identities of Cartan and Bianchi type are presented. \
\ \ \ \bigskip\newline
\textbf{2000 Mathematics Subject Classification:} 00A69, 58A15,
58B34.\bigskip\newline
\ \ \ \textbf{Keywords:} vector bundle, (generalized) Lie algebroid,
interior differential system, exterior differential calculus, exterior
differential system, Cartan identities, Bianchi identities.
\end{abstract}

\tableofcontents

\ \ 


%



\section{Introduction}

Using the notion of generalized Lie algebroid introduced in the paper $[1]$
we present the category of generalized Lie algebroids. In the
framework of Lie algebroids (see $\left[ 2\right] $) we know the following

\textbf{Theorem} (of Maurer-Cartan type) \emph{If }$((F,\nu
,N),[,]_{F},(\rho ,Id_{N}))$ \emph{is a Lie algebroid and} $d^{F}$\break
\emph{is the ex\-te\-rior differentiation operator of the exterior
differential} $\mathcal{F}(N)$\emph{-algebra}\break $(\Lambda (F,\nu
,N),+,\cdot ,\wedge ),$ \emph{then we obtain the structure equations of
Maurer-Cartan type }%
\begin{equation*}
\begin{array}{c}
d^{F}t^{\alpha }=-\displaystyle\frac{1}{2}L_{\beta \gamma }^{\alpha
}t^{\beta }\wedge t^{\gamma },~\alpha \in \overline{1,p}%
\end{array}%
\leqno(\mathcal{C}_{1})
\end{equation*}%
\emph{and\ }%
\begin{equation*}
\begin{array}{c}
d^{F}x^{i}=\rho _{\alpha }^{i}t^{\alpha },~i\in \overline{1,n},%
\end{array}%
\leqno(\mathcal{C}_{2})
\end{equation*}%
\emph{where }$\left\{ t^{\alpha },\alpha \in \overline{1,p}\right\} ~$\emph{%
is the coframe of the vector bundle }$\left( F,\nu ,N\right) .$\bigskip
\noindent

These equations are called \emph{the structure equations of Maurer-Cartan
type associa\-ted to the Lie algebroid }$\left( \left( F,\nu ,N\right) ,%
\left[ ,\right] _{F},\left( \rho ,Id_{N}\right) \right) .$

In this paper we present an exterior differential calculus for generalized
Lie algebroids. In particular, we obtain similar results with the classical
results for Lie algebroids. (see: $[4,8,9]$) A new result of Maurer-Cartan
type is presented.

We know (see $\left[ 2\right] $) that an \emph{interior differential system
(IDS)} of an Lie algebroid $\left( \left( F,\nu ,N\right) ,\left[ ,\right]
_{F},\left( \rho ,Id_{N}\right) \right) $ is a vector subbundle $\left(
E,\pi ,N\right) $ of the vector bundle $\left( F,\nu ,N\right) $. The \emph{%
IDS} $\left( E,\pi ,N\right) $ is called \emph{involutive} if $\left[ S,T%
\right] _{F}\in \Gamma \left( E,\pi ,N\right) ,~$for any $S,T\in \Gamma
\left( E,\pi ,N\right) .$ If $\left( E,\pi ,M\right) $ is an \emph{IDS} of
the Lie algebroid
\begin{equation*}
\left( \left( F,\nu ,N\right) ,\left[ ,\right] _{F},\left( \rho
,Id_{N}\right) \right) ,
\end{equation*}%
then we obtain a vector subbundle $\left( E^{0},\pi ^{0},N\right) $ of the
dual vector bundle $\left( \overset{\ast }{F},\overset{\ast }{\nu },N\right)
$ so that
\begin{equation*}
\Gamma \left( E^{0},\pi ^{0},N\right) \overset{put}{=}\left\{ \Omega \in
\Gamma \left( \overset{\ast }{F},\overset{\ast }{\nu },N\right) :\Omega
\left( S\right) =0,~\forall S\in \Gamma \left( E,\pi ,N\right) \right\} .
\end{equation*}

The vector subbundle $\left( E^{0},\pi ^{0},N\right) $ is called \emph{the
annihilator vector subbundle of the IDS }$\left( E,\pi ,N\right) .$

A characterisation of the involutivity of an \emph{IDS} (see $\left[ 2\right]
$) is presented in the following

\textbf{Theorem} (of Cartan type)\textbf{\ }\emph{Let }$\left( E,\pi
,N\right) $\emph{\ be an IDS of the Lie algebroid }%
\begin{equation*}
\left( \left( F,\nu ,N\right) ,\left[ ,\right] _{F},\left( \rho
,Id_{N}\right) \right) .
\end{equation*}

\emph{If }$\left\{ \Theta ^{r+1},...,\Theta ^{p}\right\} $\emph{\ is a base
of the }$\mathcal{F}\left( N\right) $\emph{-submodule }$\left( \Gamma \left(
E^{0},\pi ^{0},N\right) ,+,\cdot \right) $\emph{, then the IDS }$\left(
E,\pi ,N\right) $\emph{\ is involutive if and only if it exists }%
\begin{equation*}
\Omega _{\beta }^{\alpha }\in \Lambda ^{1}\left( F,\nu ,N\right) ,~\alpha
,\beta \in \overline{r+1,p}
\end{equation*}%
\emph{so that}
\begin{equation*}
d^{F}\Theta ^{\alpha }=\Sigma _{\beta \in \overline{r+1,p}}\Omega _{\beta
}^{\alpha }\wedge \Theta ^{\beta }\in \mathcal{I}\left( \Gamma \left(
E^{0},\pi ^{0},N\right) \right) .
\end{equation*}

In this paper we extend the notion of \emph{IDS} for generalized Lie
algebroids and we characterized the involutivity of an \emph{IDS} in a new theorem
of Cartan type.

The classical notion of \emph{exterior differential system (EDS)} was
studied in many papers. (see: $[3,5,6,7]$) A new point of view in the
framework of Lie algebroids is presented in the paper $\left[ 2\right] .$

Any ideal $\left( \mathcal{I},+,\cdot \right) $ of the exterior differential
algebra of the Lie algebroid
\begin{equation*}
\left( \left( F,\nu ,N\right) ,\left[ ,\right] _{F},\left( \rho
,Id_{M}\right) \right)
\end{equation*}
closed under differentiation operator $d^{F}$ $,$ namely $d^{F}\mathcal{%
I\subseteq I},$ is called \emph{differential ideal}$.$ If $\left( \mathcal{I}%
,+,\cdot \right) $ is a differential ideal of the Lie algebroid $\left(
\left( F,\nu ,N\right) ,\left[ ,\right] _{F},\left( \rho ,Id_{N}\right)
\right) $ so that it exists an \emph{IDS} $\left( E,\pi ,N\right) $ so that
for all $k\in \mathbb{N}^{\ast }$ and $\omega \in \mathcal{I}\cap \Lambda
^{k}\left( F,\nu ,N\right) $ we have $\omega \left( u_{1},...,u_{k}\right)
=0,$ for any $u_{1},...,u_{k}\in \Gamma \left( E,\pi ,N\right) ,$ then we
say that $\left( \mathcal{I},+,\cdot \right) $\emph{\ is an exterior
differential system (EDS) of the Lie algebroid }$\left( \left( F,\nu
,N\right) ,\left[ ,\right] _{F},\left( \rho ,Id_{N}\right) \right) .$

In the paper $\left[ 2\right] $ is presented the following

\textbf{Theorem} (of Cartan type) \emph{The IDS }$\left( E,\pi ,N\right) $%
\emph{\ of the Lie algebroid }%
\begin{equation*}
\left( \left( F,\nu ,N\right) ,\left[ ,\right] _{F},\left( \rho
,Id_{N}\right) \right)
\end{equation*}%
\emph{\ is involutive, if and only if the ideal generated by the }$\mathcal{F%
}\left( N\right) $\emph{-submodule }$\left( \Gamma \left( E^{0},\pi
^{0},N\right) ,+,\cdot \right) $\emph{\ is an EDS of the Lie algebroid }$%
\left( \left( F,\nu ,N\right) ,\left[ ,\right] _{F},\left( \rho
,Id_{N}\right) \right) .$

In this paper we extend the notion of \emph{EDS} to generalized Lie
algebroids. The involutivity of an \emph{IDS} in a theorem of Cartan type is
characterized. Finally, using the theory of linear connections of Ehresmann
type presented in the paper $\left[ 1\right] $, the identities of Cartan and
Bianchi type emphasize the utility of the exterior differential calculus for generalized Lie algebroids. 

\section{The category of generalized Lie algebroids}

In general, if $\mathcal{C}$ is a category, then we denote $\left\vert
\mathcal{C}\right\vert $ the class of objects and for any $A,B{\in }%
\left\vert \mathcal{C}\right\vert $, we denote $\mathcal{C}\left( A,B\right)
$ the set of morphisms of $A$ source and $B$ target. Let$\mathbf{~Vect},$ $%
\mathbf{Liealg},~\mathbf{Mod}$\textbf{,} $\mathbf{Man}$ and $\mathbf{B}^{%
\mathbf{v}}$ be the category of real vector spaces, Lie algebras, modules,
manifolds and vector bundles respectively.

We know that if $\left( E,\pi ,M\right) \in \left\vert \mathbf{B}^{\mathbf{v}%
}\right\vert ,$ $\Gamma \left( E,\pi ,M\right) =\left\{ u\in \mathbf{Man}%
\left( M,E\right) :u\circ \pi =Id_{M}\right\} $ and $\mathcal{F}\left(
M\right) =\mathbf{Man}\left( M,\mathbb{R}\right) ,$ then $\left( \Gamma
\left( E,\pi ,M\right) ,+,\cdot \right) $ is a $\mathcal{F}\left( M\right) $%
-module. If \ $\left( \varphi ,\varphi _{0}\right) \in \mathbf{B}^{\mathbf{v}%
}\left( \left( E,\pi ,M\right) ,\left( E^{\prime },\pi ^{\prime },M^{\prime
}\right) \right) $ such that $\varphi _{0}\in Iso_{\mathbf{Man}}\left(
M,M^{\prime }\right) ,$ then, using the operation
\begin{equation*}
\begin{array}{ccc}
\mathcal{F}\left( M\right) \times \Gamma \left( E^{\prime },\pi ^{\prime
},M^{\prime }\right) & ^{\underrightarrow{~\ \ \cdot ~\ \ }} & \Gamma \left(
E^{\prime },\pi ^{\prime },M^{\prime }\right) \\
\left( f,u^{\prime }\right) & \longmapsto & f\circ \varphi _{0}^{-1}\cdot
u^{\prime }%
\end{array}%
\end{equation*}%
it results that $\left( \Gamma \left( E^{\prime },\pi ^{\prime },M^{\prime
}\right) ,+,\cdot \right) $ is a $\mathcal{F}\left( M\right) $-module and we
obtain the $\mathbf{Mod}$-morphism%
\begin{equation*}
\begin{array}{ccc}
\Gamma \left( E,\pi ,M\right) & ^{\underrightarrow{~\ \ \Gamma \left(
\varphi ,\varphi _{0}\right) ~\ \ }} & \Gamma \left( E^{\prime },\pi
^{\prime },M^{\prime }\right) \\
u & \longmapsto & \Gamma \left( \varphi ,\varphi _{0}\right) u%
\end{array}%
\end{equation*}%
defined by
\begin{equation*}
\begin{array}{c}
\Gamma \left( \varphi ,\varphi _{0}\right) u\left( y\right) =\varphi \left(
u_{\varphi _{0}^{-1}\left( y\right) }\right) ,%
\end{array}%
\end{equation*}%
for any $y\in M^{\prime }.$

We know that a Lie algebroid is a vector bundle $\left( F,\nu ,N\right) \in
\left\vert \mathbf{B}^{\mathbf{v}}\right\vert $ such that there exists
\begin{equation*}
\begin{array}{c}
\left( \rho ,Id_{N}\right) \in \mathbf{B}^{\mathbf{v}}\left( \left( F,\nu
,N\right) ,\left( TN,\tau _{N},N\right) \right)%
\end{array}%
\end{equation*}%
and an operation
\begin{equation*}
\begin{array}{ccc}
\Gamma \left( F,\nu ,N\right) \times \Gamma \left( F,\nu ,N\right) & ^{%
\underrightarrow{\,\left[ ,\right] _{F}\,}} & \Gamma \left( F,\nu ,N\right)
\\
\left( u,v\right) & \longmapsto & \left[ u,v\right] _{F}%
\end{array}%
\end{equation*}%
with the following properties:

\begin{itemize}
\item[$LA_{1}$.] the equality holds good
\begin{equation*}
\begin{array}{c}
\left[ u,f\cdot v\right] _{F}=f\left[ u,v\right] _{F}+\Gamma \left( \rho
,Id_{N}\right) \left( u\right) f\cdot v,%
\end{array}%
\end{equation*}%
for all $u,v\in \Gamma \left( F,\nu ,N\right) $ and $f\in \mathcal{F}\left(
N\right) ,$

\item[$LA_{2}$.] the $4$-tuple $\left( \Gamma \left( F,\nu ,N\right)
,+,\cdot ,\left[ ,\right] _{F}\right) $ is a Lie $\mathcal{F}\left( N\right)
$-algebra$,$

\item[$LA_{3}$.] the $\mathbf{Mod}$-morphism $\Gamma \left( \rho
,Id_{N}\right) $ is a $\mathbf{LieAlg}$-morphism of $\left( \Gamma \left(
F,\nu ,N\right) ,+,\cdot ,\left[ ,\right] _{F}\right) $ source and $\left(
\Gamma \left( TN,\tau _{N},N\right) ,+,\cdot ,\left[ ,\right] _{TN}\right) $
target.
\end{itemize}

\textbf{Definition 2.1 }Let $M,N\in \left\vert \mathbf{Man}\right\vert ,$ $%
h\in Iso_{\mathbf{Man}}\left( M,N\right) $ and $\eta \in Iso_{\mathbf{Man}%
}\left( N,M\right) $.

If $\left( F,\nu ,N\right) \in \left\vert \mathbf{B}^{\mathbf{v}}\right\vert
$ so that there exists
\begin{equation*}
\begin{array}{c}
\left( \rho ,\eta \right) \in \mathbf{B}^{\mathbf{v}}\left( \left( F,\nu
,N\right) ,\left( TM,\tau _{M},M\right) \right)%
\end{array}%
\end{equation*}%
and an operation
\begin{equation*}
\begin{array}{ccc}
\Gamma \left( F,\nu ,N\right) \times \Gamma \left( F,\nu ,N\right) & ^{%
\underrightarrow{\left[ ,\right] _{F,h}}} & \Gamma \left( F,\nu ,N\right) \\
\left( u,v\right) & \longmapsto & \left[ u,v\right] _{F,h}%
\end{array}%
\end{equation*}%
with the following properties:\bigskip

\noindent $\qquad GLA_{1}$. the equality holds good
\begin{equation*}
\begin{array}{c}
\left[ u,f\cdot v\right] _{F,h}=f\left[ u,v\right] _{F,h}+\Gamma \left(
Th\circ \rho ,h\circ \eta \right) \left( u\right) f\cdot v,%
\end{array}%
\end{equation*}%
\qquad \quad\ \ for all $u,v\in \Gamma \left( F,\nu ,N\right) $ and $f\in
\mathcal{F}\left( N\right) .$

\medskip $GLA_{2}$. the $4$-tuple $\left( \Gamma \left( F,\nu ,N\right)
,+,\cdot ,\left[ ,\right] _{F,h}\right) $ is a Lie $\mathcal{F}\left(
N\right) $-algebra,

$GLA_{3}$. the $\mathbf{Mod}$-morphism $\Gamma \left( Th\circ \rho ,h\circ
\eta \right) $ is a $\mathbf{LieAlg}$-morphism of
\begin{equation*}
\left( \Gamma \left( F,\nu ,N\right) ,+,\cdot ,\left[ ,\right] _{F,h}\right)
\end{equation*}
source and
\begin{equation*}
\left( \Gamma \left( TN,\tau _{N},N\right) ,+,\cdot ,\left[ ,\right]
_{TN}\right)
\end{equation*}
target, \medskip \noindent then we will say that \emph{the triple }$\left(
\left( F,\nu ,N\right) ,\left[ ,\right] _{F,h},\left( \rho ,\eta \right)
\right) $ \emph{is a generalized Lie algebroid. }The couple $\left( \left[ ,%
\right] _{F,h},\left( \rho ,\eta \right) \right) $ will be called \emph{%
generalized Lie algebroid structure.}

\textbf{Definition 2.2 }We define the set of morphisms of
\begin{equation*}
\left( \left( F,\nu ,N\right) ,\left[ ,\right] _{F,h},\left( \rho ,\eta
\right) \right)
\end{equation*}%
source and
\begin{equation*}
\left( \left( F^{\prime },\nu ^{\prime },N^{\prime }\right) ,\left[ ,\right]
_{F^{\prime },h^{\prime }},\left( \rho ^{\prime },\eta ^{\prime }\right)
\right)
\end{equation*}%
target as being the set
\begin{equation*}
\begin{array}{c}
\left\{ \left( \varphi ,\varphi _{0}\right) \in \mathbf{B}^{\mathbf{v}%
}\left( \left( F,\nu ,N\right) ,\left( F^{\prime },\nu ^{\prime },N^{\prime
}\right) \right) \right\}%
\end{array}%
\end{equation*}%
such that $\varphi _{0}\in Iso_{\mathbf{Man}}\left( N,N^{\prime }\right) $
and the $\mathbf{Mod}$-morphism $\Gamma \left( \varphi ,\varphi _{0}\right) $
is a $\mathbf{LieAlg}$-morphism of
\begin{equation*}
\left( \Gamma \left( F,\nu ,N\right) ,+,\cdot ,\left[ ,\right] _{F,h}\right)
\end{equation*}%
source and
\begin{equation*}
\left( \Gamma \left( F^{\prime },\nu ^{\prime },N^{\prime }\right) ,+,\cdot
, \left[ ,\right] _{F^{\prime },h^{\prime }}\right)
\end{equation*}%
target.

So, we can discuss about \emph{the category }$\mathbf{GLA}$\emph{\ of
generalized Lie algebroids.} Examples of objects of this category are
presented in the paper $\left[ 1\right] .$ We remark that $\mathbf{GLA}$ is
a subcategory of the category $\mathbf{B}^{\mathbf{v}}.$

Let $\left( \left( F,\nu ,N\right) ,\left[ ,\right] _{F,h},\left( \rho ,\eta
\right) \right) $ be an arbitrary object of the category $\mathbf{GLA}$.

\begin{itemize}
\item Locally, for any $\alpha ,\beta \in \overline{1,p},$ we set $\left[
t_{\alpha },t_{\beta }\right] _{F,h}=L_{\alpha \beta }^{\gamma }t_{\gamma }.$
We easily obtain that $L_{\alpha \beta }^{\gamma }=-L_{\beta \alpha
}^{\gamma },~$for any $\alpha ,\beta ,\gamma \in \overline{1,p}.$
\end{itemize}

The real local functions $L_{\alpha \beta }^{\gamma },~\alpha ,\beta ,\gamma
\in \overline{1,p}$ will be called the \emph{structure functions of the
generalized Lie algebroid }$\left( \left( F,\nu ,N\right) ,\left[ ,\right]
_{F,h},\left( \rho ,\eta \right) \right) .$

\begin{itemize}
\item We assume the following diagrams:%
\begin{equation*}
\begin{array}[b]{ccccc}
F & ^{\underrightarrow{~\ \ \ \rho ~\ \ }} & TM & ^{\underrightarrow{~\ \ \
Th~\ \ }} & TN \\
~\downarrow \nu &  & ~\ \ \ \downarrow \tau _{M} &  & ~\ \ \ \downarrow \tau
_{N} \\
N & ^{\underrightarrow{~\ \ \ \eta ~\ \ }} & M & ^{\underrightarrow{~\ \ \
h~\ \ }} & N \\
&  &  &  &  \\
\left( \chi ^{\tilde{\imath}},z^{\alpha }\right) &  & \left(
x^{i},y^{i}\right) &  & \left( \chi ^{\tilde{\imath}},z^{\tilde{\imath}%
}\right)%
\end{array}%
\end{equation*}

where $i,\tilde{\imath}\in \overline{1,m}$ and $\alpha \in \overline{1,p}.$

If%
\begin{equation*}
\left( \chi ^{\tilde{\imath}},z^{\alpha }\right) \longrightarrow \left( \chi
^{\tilde{\imath}\prime }\left( \chi ^{\tilde{\imath}}\right) ,z^{\alpha
\prime }\left( \chi ^{\tilde{\imath}},z^{\alpha }\right) \right) ,
\end{equation*}%
\begin{equation*}
\left( x^{i},y^{i}\right) \longrightarrow \left( x^{i%
{\acute{}}%
}\left( x^{i}\right) ,y^{i%
{\acute{}}%
}\left( x^{i},y^{i}\right) \right)
\end{equation*}%
and
\begin{equation*}
\left( \chi ^{\tilde{\imath}},z^{\tilde{\imath}}\right) \longrightarrow
\left( \chi ^{\tilde{\imath}\prime }\left( \chi ^{\tilde{\imath}}\right) ,z^{%
\tilde{\imath}\prime }\left( \chi ^{\tilde{\imath}},z^{\tilde{\imath}%
}\right) \right) ,
\end{equation*}%
then
\begin{equation*}
\begin{array}[b]{c}
z^{\alpha
{\acute{}}%
}=\Lambda _{\alpha }^{\alpha
{\acute{}}%
}z^{\alpha }%
\end{array}%
,
\end{equation*}%
\begin{equation*}
\begin{array}[b]{c}
y^{i%
{\acute{}}%
}=\frac{\partial x^{i%
{\acute{}}%
}}{\partial x^{i}}y^{i}%
\end{array}%
\end{equation*}%
and
\begin{equation*}
\begin{array}{c}
z^{\tilde{\imath}\prime }=\frac{\partial \chi ^{\tilde{\imath}\prime }}{%
\partial \chi ^{\tilde{\imath}}}z^{\tilde{\imath}}.%
\end{array}%
\end{equation*}

\item We assume that $\left( \theta ,\mu \right) \overset{put}{=}\left(
Th\circ \rho ,h\circ \eta \right) $. If $z^{\alpha }t_{\alpha }\in \Gamma
\left( F,\nu ,N\right) $ is arbitrary, then
\begin{equation*}
\begin{array}[t]{l}
\displaystyle%
\begin{array}{c}
\Gamma \left( Th\circ \rho ,h\circ \eta \right) \left( z^{\alpha }t_{\alpha
}\right) f\left( h\circ \eta \left( \varkappa \right) \right) =\vspace*{1mm}
\\
=\left( \theta _{\alpha }^{\tilde{\imath}}z^{\alpha }\frac{\partial f}{%
\partial \varkappa ^{\tilde{\imath}}}\right) \left( h\circ \eta \left(
\varkappa \right) \right) =\left( \left( \rho _{\alpha }^{i}\circ h\right)
\left( z^{\alpha }\circ h\right) \frac{\partial f\circ h}{\partial x^{i}}%
\right) \left( \eta \left( \varkappa \right) \right) ,%
\end{array}%
\end{array}%
\leqno(2.1)
\end{equation*}%
for any $f\in \mathcal{F}\left( N\right) $ and $\varkappa \in N.$
\end{itemize}

The coefficients $\rho _{\alpha }^{i}$ respectively $\theta _{\alpha }^{%
\tilde{\imath}}$ change to $\rho _{\alpha
{\acute{}}%
}^{i%
{\acute{}}%
}$ respectively $\theta _{\alpha
{\acute{}}%
}^{\tilde{\imath}%
{\acute{}}%
}$ according to the rule:
\begin{equation*}
\begin{array}{c}
\rho _{\alpha
{\acute{}}%
}^{i%
{\acute{}}%
}=\Lambda _{\alpha
{\acute{}}%
}^{\alpha }\rho _{\alpha }^{i}\displaystyle\frac{\partial x^{i%
{\acute{}}%
}}{\partial x^{i}},%
\end{array}%
\leqno(2.2)
\end{equation*}%
respectively%
\begin{equation*}
\begin{array}{c}
\theta _{\alpha
{\acute{}}%
}^{\tilde{\imath}%
{\acute{}}%
}=\Lambda _{\alpha
{\acute{}}%
}^{\alpha }\theta _{\alpha }^{\tilde{\imath}}\displaystyle\frac{\partial
\varkappa ^{\tilde{\imath}%
{\acute{}}%
}}{\partial \varkappa ^{\tilde{\imath}}},%
\end{array}%
\leqno(2.3)
\end{equation*}%
where
\begin{equation*}
\left\Vert \Lambda _{\alpha
{\acute{}}%
}^{\alpha }\right\Vert =\left\Vert \Lambda _{\alpha }^{\alpha
{\acute{}}%
}\right\Vert ^{-1}.
\end{equation*}

\textit{Remark 2.1 }\emph{The following equalities hold good:}%
\begin{equation*}
\begin{array}{c}
\displaystyle\rho _{\alpha }^{i}\circ h\frac{\partial f\circ h}{\partial
x^{i}}=\left( \theta _{\alpha }^{\tilde{\imath}}\frac{\partial f}{\partial
\varkappa ^{\tilde{\imath}}}\right) \circ h,\forall f\in \mathcal{F}\left(
N\right) .%
\end{array}%
\leqno(2.4)
\end{equation*}%
\emph{and }%
\begin{equation*}
\begin{array}{c}
\displaystyle\left( L_{\alpha \beta }^{\gamma }\circ h\right) \left( \rho
_{\gamma }^{k}\circ h\right) =\left( \rho _{\alpha }^{i}\circ h\right) \frac{%
\partial \left( \rho _{\beta }^{k}\circ h\right) }{\partial x^{i}}-\left(
\rho _{\beta }^{j}\circ h\right) \frac{\partial \left( \rho _{\alpha
}^{k}\circ h\right) }{\partial x^{j}}.%
\end{array}%
\leqno(2.5)
\end{equation*}

\section{Interior Differential Systems}

Let $\left( \left( F,\nu ,N\right) ,\left[ ,\right] _{F,h},\left( \rho ,\eta
\right) \right) $ be an object of the category $\mathbf{GLA}$.

Let $\mathcal{AF}_{F}$ be a vector fibred $\left( n+p\right) $-atlas for the
vector bundle $\left( F,\nu ,N\right) $ and let $\mathcal{AF}_{TM}$ be a
vector fibred $\left( m+m\right) $-atlas for the vector bundle $\left(
TM,\tau _{M},M\right) $.

Let $\left( h^{\ast }F,h^{\ast }\nu ,M\right) $ be the pull-back vector
bundle through $h.$

If $\left( U,\xi _{U}\right) \in \mathcal{AF}_{TM}$ and $\left(
V,s_{V}\right) \in \mathcal{AF}_{F}$ such that $U\cap h^{-1}\left( V\right)
\neq \phi $, then we define the application%
\begin{equation*}
\begin{array}{ccc}
h^{\ast }\nu ^{-1}(U{\cap }h^{-1}(V))) & {}^{\underrightarrow{\bar{s}_{U{%
\cap }h^{-1}(V)}}} & \left( U{\cap }h^{-1}(V)\right) {\times }\mathbb{R}^{p}
\\
\left( \varkappa ,z\left( h\left( \varkappa \right) \right) \right) &
\longmapsto & \left( \varkappa ,t_{V,h\left( \varkappa \right) }^{-1}z\left(
h\left( \varkappa \right) \right) \right).%
\end{array}%
\end{equation*}

\textbf{Proposition 3.1 }\emph{The set}%
\begin{equation*}
\begin{array}{c}
\overline{\mathcal{AF}}_{F}\overset{put}{=}\underset{U\cap h^{-1}\left(
V\right) \neq \phi }{\underset{\left( U,\xi _{U}\right) \in \mathcal{AF}%
_{TM},~\left( V,s_{V}\right) \in \mathcal{AF}_{F}}{\tbigcup }}\left\{ \left(
U\cap h^{-1}\left( V\right) ,\bar{s}_{U{\cap }h^{-1}(V)}\right) \right\}%
\end{array}%
\end{equation*}%
\emph{is a vector fibred }$m+p$\emph{-atlas for the vector bundle }$\left(
h^{\ast }F,h^{\ast }\nu ,M\right) .$

\emph{If }$z=z^{\alpha }t_{\alpha }\in \Gamma \left( F,\nu ,N\right) ,$
\emph{then\ we obtain the section }%
\begin{equation*}
\begin{array}{c}
Z=\left( z^{\alpha }\circ h\right) T_{\alpha }\in \Gamma \left( h^{\ast
}F,h^{\ast }\nu ,M\right)%
\end{array}%
\end{equation*}%
\emph{such that }$Z\left( x\right) =z\left( h\left( x\right) \right) ,$
\emph{for any }$x\in U\cap h^{-1}\left( V\right) .$\bigskip

\textbf{Theorem 3.1 }\emph{Let} $\Big({\overset{h^{\ast }F}{\rho }},Id_{M}%
\Big)$ \emph{be the }$\mathbf{B}^{\mathbf{v}}$\emph{-morphism of }$\left(
h^{\ast }F,h^{\ast }\nu ,M\right) $\ \emph{source and} $\left( TM,\tau
_{M},M\right) $\ \emph{target, where}%
\begin{equation*}
\begin{array}{rcl}
h^{\ast }F & ^{\underrightarrow{\overset{h^{\ast }F}{\rho }}} & TM \\
\displaystyle Z^{\alpha }T_{\alpha }\left( x\right) & \longmapsto & %
\displaystyle\left( Z^{\alpha }\cdot \rho _{\alpha }^{i}\circ h\right) \frac{%
\partial }{\partial x^{i}}\left( x\right)%
\end{array}%
\leqno(3.1)
\end{equation*}

\emph{Using the operation}
\begin{equation*}
\begin{array}{ccc}
\Gamma \left( h^{\ast }F,h^{\ast }\nu ,M\right) \times \Gamma \left( h^{\ast
}F,h^{\ast }\nu ,M\right) & ^{\underrightarrow{~\ \ \left[ ,\right]
_{h^{\ast }F}~\ \ }} & \Gamma \left( h^{\ast }F,h^{\ast }\nu ,M\right)%
\end{array}%
\end{equation*}%
\emph{defined by}%
\begin{equation*}
\begin{array}{ll}
\left[ T_{\alpha },T_{\beta }\right] _{h^{\ast }F} & =\left( L_{\alpha \beta
}^{\gamma }\circ h\right) T_{\gamma },\vspace*{1mm} \\
\left[ T_{\alpha },fT_{\beta }\right] _{h^{\ast }F} & \displaystyle=f\left(
L_{\alpha \beta }^{\gamma }\circ h\right) T_{\gamma }+\left( \rho _{\alpha
}^{i}\circ h\right) \frac{\partial f}{\partial x^{i}}T_{\beta },\vspace*{1mm}
\\
\left[ fT_{\alpha },T_{\beta }\right] _{h^{\ast }F} & =-\left[ T_{\beta
},fT_{\alpha }\right] _{h^{\ast }F},%
\end{array}%
\leqno(3.2)
\end{equation*}%
\emph{for any} $f\in \mathcal{F}\left( M\right) ,$ \emph{it results that}
\begin{equation*}
\begin{array}{c}
\left( \left( h^{\ast }F,h^{\ast }\nu ,M\right) ,\left[ ,\right] _{h^{\ast
}F},\left( \overset{h^{\ast }F}{\rho },Id_{M}\right) \right)%
\end{array}%
\end{equation*}%
is a Lie algebroid which is called \emph{the pull-back Lie algebroid of the
generalized Lie algebroid }$\left( \left( F,\nu ,N\right) ,\left[ ,\right]
_{F,h},\left( \rho ,\eta \right) \right) .$

\textbf{Definition 3.1 }Any vector subbundle $\left( E,\pi ,M\right) $ of
the pull-back vector bundle $\left( h^{\ast }F,h^{\ast }\nu ,M\right) $ will
be called \emph{interior differential system (IDS) of the generalized Lie
algebroid }%
\begin{equation*}
\left( \left( F,\nu ,N\right) ,\left[ ,\right] _{F,h},\left( \rho ,\eta
\right) \right) .
\end{equation*}

In particular, if $h=Id_{N}=\eta $, then we obtain the definition of the IDS
of a Lie algebroid. (see $\left[ {2}\right] $)

\emph{Remark 3.1 }If $\left( E,\pi ,M\right) $ is an IDS of the generalized
Lie algebroid
\begin{equation*}
\left( \left( F,\nu ,N\right) ,\left[ ,\right] _{F,h},\left( \rho ,\eta
\right) \right) ,
\end{equation*}%
then we obtain a vector subbundle $\left( E^{0},\pi ^{0},M\right) $ of the
vector bundle $\left( \overset{\ast }{h^{\ast }F},\overset{\ast }{h^{\ast
}\nu },M\right) $ such that
\begin{equation*}
\Gamma \left( E^{0},\pi ^{0},M\right) \overset{put}{=}\left\{ \Omega \in
\Gamma \left( \overset{\ast }{h^{\ast }F},\overset{\ast }{h^{\ast }\nu }%
,M\right) :\Omega \left( S\right) =0,~\forall S\in \Gamma \left( E,\pi
,M\right) \right\} .
\end{equation*}

The vector subbundle $\left( E^{0},\pi ^{0},M\right) $ will be called \emph{%
the annihilator vector subbundle of the IDS }$\left( E,\pi ,M\right) .$

\textbf{Proposition 3.1 }\emph{If }$\left( E,\pi ,M\right) $\emph{\ is an
IDS of the generalized Lie algebroid }%
\begin{equation*}
\left( \left( F,\nu ,N\right) ,\left[ ,\right] _{F,h},\left( \rho ,\eta
\right) \right)
\end{equation*}%
\emph{such that }$\Gamma \left( E,\pi ,M\right) =\left\langle
S_{1},...,S_{r}\right\rangle $\emph{, then it exists }$\Theta
^{r+1},...,\Theta ^{p}\in \Gamma \left( \overset{\ast }{h^{\ast }F},\overset{%
\ast }{h^{\ast }\nu },M\right) $\emph{\ linearly independent such that }$%
\Gamma \left( E^{0},\pi ^{0},M\right) =\left\langle \Theta ^{r+1},...,\Theta
^{p}\right\rangle .$

\textbf{Definition 3.2 }The IDS $\left( E,\pi ,M\right) $ of the generalized
Lie algebroid
\begin{equation*}
\left( \left( F,\nu ,N\right) ,\left[ ,\right] _{F,h},\left( \rho ,\eta
\right) \right)
\end{equation*}%
will be called \emph{involutive} if $\left[ S,T\right] _{h^{\ast }F}\in
\Gamma \left( E,\pi ,M\right) ,~$for any $S,T\in \Gamma \left( E,\pi
,M\right) .$

\textbf{Proposition 3.2 }\emph{If }$\left( E,\pi ,M\right) $\emph{\ is an
IDS of the generalized Lie algebroid }%
\begin{equation*}
\left( \left( F,\nu ,N\right) ,\left[ ,\right] _{F,h},\left( \rho ,\eta
\right) \right)
\end{equation*}%
\emph{and }$\left\{ S_{1},...,S_{r}\right\} $\emph{\ is a base for the }$%
\mathcal{F}\left( M\right) $\emph{-submodule }$\left( \Gamma \left( E,\pi
,M\right) ,+,\cdot \right) $\emph{\ then }$\left( E,\pi ,M\right) $\emph{\
is involutive if and only if }$\left[ S_{a},S_{b}\right] _{h^{\ast }F}\in
\Gamma \left( E,\pi ,M\right) ,~$for any $a,b\in \overline{1,r}.$

\section{Exterior differential calculus}

We propose an exterior differential calculus in the general framework of
generalized Lie algebroids. As any Lie algebroid can be regarded as a
generalized Lie algebroid, in particular, we obtain a new point of view over
the exterior differential calculus for Lie algebroids. Let $\left( \left(
F,\nu ,N\right) ,\left[ ,\right] _{F,h},\left( \rho ,\eta \right) \right)
\in \left\vert \mathbf{GLA}\right\vert $ be.

\textbf{Definition 4.1 }For any $q\in \mathbb{N}$ we denote by $\left(
\Sigma _{q},\circ \right) $ the permutations group of the set $\left\{
1,2,...,q\right\} .$

\textbf{Definition 4.2 }We denoted by $\Lambda ^{q}\left( F,\nu ,N\right) $
the set of $q$-linear applications
\begin{equation*}
\begin{array}{ccc}
\Gamma \left( F,\nu ,N\right) ^{q} & ^{\underrightarrow{\ \ \omega \ \ }} &
\mathcal{F}\left( N\right) \\
\left( z_{1},...,z_{q}\right) & \longmapsto & \omega \left(
z_{1},...,z_{q}\right)%
\end{array}%
\end{equation*}%
such that
\begin{equation*}
\begin{array}{c}
\omega \left( z_{\sigma \left( 1\right) },...,z_{\sigma \left( q\right)
}\right) =sgn\left( \sigma \right) \cdot \omega \left( z_{1},...,z_{q}\right)%
\end{array}%
\end{equation*}%
for any $z_{1},...,z_{q}\in \Gamma \left( F,\nu ,N\right) $ and for any $%
\sigma \in \Sigma _{q}$.

The elements of $\Lambda ^{q}\left( F,\nu ,N\right) $ will be called \emph{%
differential forms of degree }$q$ or \emph{differential }$q$\emph{-forms}$.$

\textit{\noindent Remark 4.1}\textbf{\ }If $\omega \in \Lambda ^{q}\left(
F,\nu ,N\right) $, then $\omega \left( z_{1},...,z,...,z,...z_{q}\right) =0.$
Therefore, if $\omega \in \Lambda ^{q}\left( F,\nu ,N\right) $, then
\begin{equation*}
\begin{array}{c}
\omega \left( z_{1},...,z_{i},...,z_{j},...z_{q}\right) =-\omega \left(
z_{1},...,z_{j},...,z_{i},...z_{q}\right) .%
\end{array}%
\end{equation*}

\textbf{Theorem 4.1 }\emph{If }$q\in N$\emph{, then }$\left( \Lambda
^{q}\left( F,\nu ,N\right) ,+,\cdot \right) $\emph{\ is a} $\mathcal{F}%
\left( N\right) $\emph{-module.}

\textbf{Definition 4.3 }If $\omega \in \Lambda ^{q}\left( F,\nu ,N\right) $
and $\theta \in \Lambda ^{r}\left( F,\nu ,N\right) $, then the $\left(
q+r\right) $-form $\omega \wedge \theta $ defined by%
\begin{equation*}
\begin{array}{ll}
\omega \wedge \theta \left( z_{1},...,z_{q+r}\right) & =\underset{\sigma
\left( q+1\right) <...<\sigma \left( q+r\right) }{\underset{\sigma \left(
1\right) <...<\sigma \left( q\right) }{\tsum }}sgn\left( \sigma \right)
\omega \left( z_{\sigma \left( 1\right) },...,z_{\sigma \left( q\right)
}\right) \theta \left( z_{\sigma \left( q+1\right) },...,z_{\sigma \left(
q+r\right) }\right) \vspace{3mm} \\
& \displaystyle=\frac{1}{q!r!}\underset{\sigma \in \Sigma _{q+r}}{\tsum }%
sgn\left( \sigma \right) \omega \left( z_{\sigma \left( 1\right)
},...,z_{\sigma \left( q\right) }\right) \theta \left( z_{\sigma \left(
q+1\right) },...,z_{\sigma \left( q+r\right) }\right) ,%
\end{array}%
\end{equation*}%
for any $z_{1},...,z_{q+r}\in \Gamma \left( F,\nu ,N\right) ,$ will be
called \emph{the exterior product of the forms }$\omega $ \emph{and}~$\theta
.$

Using the previous definition, we obtain

\textbf{Theorem 4.2 }\emph{The following affirmations hold good:} \medskip

\noindent 1. \emph{If }$\omega \in \Lambda ^{q}\left( F,\nu ,N\right) $\emph{%
\ and }$\theta \in \Lambda ^{r}\left( F,\nu ,N\right) $\emph{, then}%
\begin{equation*}
\begin{array}{c}
\omega \wedge \theta =\left( -1\right) ^{q\cdot r}\theta \wedge \omega .%
\end{array}%
\leqno(4.1)
\end{equation*}

\noindent 2. \emph{For any }$\omega \in \Lambda ^{q}\left( F,\nu ,N\right) $%
\emph{, }$\theta \in \Lambda ^{r}\left( F,\nu ,N\right) $\emph{\ and }$\eta
\in \Lambda ^{s}\left( F,\nu ,N\right) $\emph{\ we obtain}%
\begin{equation*}
\begin{array}{c}
\left( \omega \wedge \theta \right) \wedge \eta =\omega \wedge \left( \theta
\wedge \eta \right) .%
\end{array}%
\leqno(4.2)
\end{equation*}

\noindent 3. \emph{For any }$\omega ,\theta \in \Lambda ^{q}\left( F,\nu
,N\right) $\emph{\ and }$\eta \in \Lambda ^{s}\left( F,\nu ,N\right) $\emph{%
\ we obtain }%
\begin{equation*}
\begin{array}{c}
\left( \omega +\theta \right) \wedge \eta =\omega \wedge \eta +\theta \wedge
\eta .%
\end{array}%
\leqno(4.3)
\end{equation*}

\noindent 4. \emph{For any }$\omega \in \Lambda ^{q}\left( F,\nu ,N\right) $%
\emph{\ and }$\theta ,\eta \in \Lambda ^{s}\left( F,\nu ,N\right) $\emph{\
we obtain }%
\begin{equation*}
\begin{array}{c}
\omega \wedge \left( \theta +\eta \right) =\omega \wedge \theta +\omega
\wedge \eta .%
\end{array}%
\leqno(4.4)
\end{equation*}

\noindent 5. \emph{For any} $f\in \mathcal{F}\left( N\right) $, $\omega \in
\Lambda ^{q}\left( F,\nu ,N\right) $ \emph{and }$\theta \in \Lambda
^{s}\left( F,\nu ,N\right) $\emph{\ we obtain }%
\begin{equation*}
\begin{array}{c}
\left( f\cdot \omega \right) \wedge \theta =f\cdot \left( \omega \wedge
\theta \right) =\omega \wedge \left( f\cdot \theta \right) .%
\end{array}%
\leqno(4.5)
\end{equation*}%
\noindent

\textbf{Theorem 4.3 }\emph{If }%
\begin{equation*}
\Lambda \left( F,\nu ,N\right) =\underset{q\geq 0}{\oplus }\Lambda
^{q}\left( F,\nu ,N\right) ,
\end{equation*}%
\emph{then }%
\begin{equation*}
\begin{array}{c}
\left( \Lambda \left( F,\nu ,N\right) ,+,\cdot ,\wedge \right)%
\end{array}%
\end{equation*}%
\emph{is a} $\mathcal{F}\left( N\right) $\emph{-algebra.} This algebra will
be called \emph{the exterior differential algebra of the vector bundle }$%
\left( F,\nu ,N\right) .$

\emph{Remark 4.2}If $\left\{ t^{\alpha },~\alpha \in \overline{1,p}\right\} $
is the coframe associated to the frame $\left\{ t_{\alpha },~\alpha \in
\overline{1,p}\right\} $ of the vector bundle $\left( F,\nu ,N\right) $ in
the vector local $\left( n+p\right) $-chart $U$, then
\begin{equation*}
\begin{array}{c}
t^{\alpha _{1}}\wedge ...\wedge t^{\alpha _{q}}\left( z_{1}^{\alpha
}t_{\alpha },...,z_{q}^{\alpha }t_{\alpha }\right) =\frac{1}{q!}\det
\left\Vert
\begin{array}{ccc}
z_{1}^{\alpha _{1}} & ... & z_{1}^{\alpha _{q}} \\
... & ... & ... \\
z_{q}^{\alpha _{1}} & ... & z_{q}^{\alpha _{q}}%
\end{array}%
\right\Vert ,%
\end{array}%
\end{equation*}%
for any $q\in \overline{1,p}.$

\emph{Remark 4.3}If $\left\{ t^{\alpha },~\alpha \in \overline{1,p}\right\} $
is the coframe associated to the frame $\left\{ t_{\alpha },~\alpha \in
\overline{1,p}\right\} $ of the vector bundle $\left( F,\nu ,N\right) $ in
the vector local $\left( n+p\right) $-chart $U$, then, for any $q\in
\overline{1,p}$ we define $C_{p}^{q}$ exterior differential forms of the
type
\begin{equation*}
\begin{array}{c}
t^{\alpha _{1}}\wedge ...\wedge t^{\alpha _{q}}%
\end{array}%
\end{equation*}%
such that $1\leq \alpha _{1}<...<\alpha _{q}\leq p.$

The set%
\begin{equation*}
\begin{array}{c}
\left\{ t^{\alpha _{1}}\wedge ...\wedge t^{\alpha _{q}},1\leq \alpha
_{1}<...<\alpha _{q}\leq p\right\}%
\end{array}%
\end{equation*}%
is a base for the $\mathcal{F}\left( N\right) $-module
\begin{equation*}
\left( \Lambda ^{q}\left( F,\nu ,N\right) ,+,\cdot \right) .
\end{equation*}

Therefore, if $\omega \in \Lambda ^{q}\left( F,\nu ,N\right) $, then%
\begin{equation*}
\begin{array}{c}
\omega =\omega _{\alpha _{1}...\alpha _{q}}t^{\alpha _{1}}\wedge ...\wedge
t^{\alpha _{q}}.%
\end{array}%
\end{equation*}

In particular, if $\omega $ is an exterior differential $p$-form $\omega ,$
then we can written
\begin{equation*}
\begin{array}{c}
\omega =a\cdot t^{1}\wedge ...\wedge t^{p},%
\end{array}%
\end{equation*}%
where $a\in \mathcal{F}\left( N\right) .$

\textbf{Definition 4.4 }If
\begin{equation*}
\begin{array}{c}
\omega =\omega _{\alpha _{1}...\alpha _{q}}t^{\alpha _{1}}\wedge ...\wedge
t^{\alpha _{q}}\in \Lambda ^{q}\left( F,\nu ,N\right)%
\end{array}%
\end{equation*}%
such that
\begin{equation*}
\begin{array}{c}
\omega _{\alpha _{1}...\alpha _{q}}\in C^{r}\left( N\right) ,%
\end{array}%
\end{equation*}%
for any $1\leq \alpha _{1}<...<\alpha _{q}\leq p$, then we will say that
\emph{the }$q$\emph{-form }$\omega $\emph{\ is differentiable of }$C^{r}$%
\emph{-class.}

\textbf{Definition 4.5 }For any $z\in \Gamma \left( F,\nu ,N\right) $, the $%
\mathcal{F}\left( N\right) $-multilinear application
\begin{equation*}
\begin{array}{c}
\begin{array}{rcl}
\Lambda \left( F,\nu ,N\right) & ^{\underrightarrow{~\ \ L_{z}~\ \ }} &
\Lambda \left( F,\nu ,N\right)%
\end{array}%
,%
\end{array}%
\end{equation*}%
defined by%
\begin{equation*}
\begin{array}{c}
L_{z}\left( f\right) =\Gamma \left( Th\circ \rho ,h\circ \eta \right)
z\left( f\right) ,~\forall f\in \mathcal{F}\left( N\right)%
\end{array}%
\end{equation*}%
and
\begin{equation*}
\begin{array}{cl}
L_{z}\omega \left( z_{1},...,z_{q}\right) & =\Gamma \left( Th\circ \rho
,h\circ \eta \right) z\left( \omega \left( \left( z_{1},...,z_{q}\right)
\right) \right) \\
& -\overset{q}{\underset{i=1}{\tsum }}\omega \left( \left( z_{1},...,\left[
z,z_{i}\right] _{F,h},...,z_{q}\right) \right) ,%
\end{array}%
\end{equation*}%
for any $\omega \in \Lambda ^{q}\mathbf{\ }\left( F,\nu ,N\right) $ and $%
z_{1},...,z_{q}\in \Gamma \left( F,\nu ,N\right) ,$ will be called \emph{the
covariant Lie derivative with respect to the section }$z.$

\textbf{Theorem 4.4 }\emph{If }$z\in \Gamma \left( F,\nu ,N\right) ,$ $%
\omega \in \Lambda ^{q}\left( F,\nu ,N\right) $\emph{\ and }$\theta \in
\Lambda ^{r}\left( F,\nu ,N\right) $\emph{, then}%
\begin{equation*}
\begin{array}{c}
L_{z}\left( \omega \wedge \theta \right) =L_{z}\omega \wedge \theta +\omega
\wedge L_{z}\theta .%
\end{array}%
\leqno(4.6)
\end{equation*}

\emph{Proof. }Let $z_{1},...,z_{q+r}\in \Gamma \left( F,\nu ,N\right) $ be
arbitrary. Since
\begin{equation*}
\begin{array}{l}
L_{z}\left( \omega \wedge \theta \right) \left( z_{1},...,z_{q+r}\right)
=\Gamma \left( Th\circ \rho ,h\circ \eta \right) z\left( \left( \omega
\wedge \theta \right) \left( z_{1},...,z_{q+r}\right) \right) \\
-\overset{q+r}{\underset{i=1}{\tsum }}\left( \omega \wedge \theta \right)
\left( \left( z_{1},...,\left[ z,z_{i}\right] _{F,h},...,z_{q+r}\right)
\right) \\
=\Gamma \left( Th\circ \rho ,h\circ \eta \right) z\left( \underset{\sigma
\left( q+1\right) <...<\sigma \left( q+r\right) }{\underset{\sigma \left(
1\right) <...<\sigma \left( q\right) }{\tsum }}sgn\left( \sigma \right)
\cdot \omega \left( z_{\sigma \left( 1\right) },...,z_{\sigma \left(
q\right) }\right) \right. \\
\qquad \left. \cdot \theta \left( z_{\sigma \left( q+1\right)
},...,z_{\sigma \left( q+r\right) }\right) \right) -\overset{q+r}{\underset{%
i=1}{\tsum }}\left( \omega \wedge \theta \right) \left( \left( z_{1},...,%
\left[ z,z_{i}\right] _{F,h},...,z_{q+r}\right) \right) \\
=\underset{\sigma \left( q+1\right) <...<\sigma \left( q+r\right) }{\underset%
{\sigma \left( 1\right) <...<\sigma \left( q\right) }{\tsum }}sgn\left(
\sigma \right) \cdot \Gamma \left( Th\circ \rho ,h\circ \eta \right) z\left(
\omega \left( z_{\sigma \left( 1\right) },...,z_{\sigma \left( q\right)
}\right) \right) \\
\qquad \cdot \theta \left( z_{\sigma \left( q+1\right) },...,z_{\sigma
\left( q+r\right) }\right) +\underset{\sigma \left( q+1\right) <...<\sigma
\left( q+r\right) }{\underset{\sigma \left( 1\right) <...<\sigma \left(
q\right) }{\tsum }}sgn\left( \sigma \right) \cdot \omega \left( z_{\sigma
\left( 1\right) },...,z_{\sigma \left( q\right) }\right) \\
\qquad \cdot \Gamma \left( Th\circ \rho ,h\circ \eta \right) z\left( \theta
\left( z_{\sigma \left( q+1\right) },...,z_{\sigma \left( q+r\right)
}\right) \right) -\underset{\sigma \left( q+1\right) <...<\sigma \left(
q+r\right) }{\underset{\sigma \left( 1\right) <...<\sigma \left( q\right) }{%
\tsum }}sgn\left( \sigma \right) \\
\qquad \cdot \overset{q}{\underset{i=1}{\tsum }}\omega \left( z_{\sigma
\left( 1\right) },...,\left[ z,z_{\sigma \left( i\right) }\right]
_{F,h},...,z_{\sigma \left( q+r\right) }\right) \cdot \theta \left(
z_{\sigma \left( q+1\right) },...,z_{\sigma \left( q+r\right) }\right)%
\end{array}%
\end{equation*}%
\begin{equation*}
\begin{array}{l}
-\underset{\sigma \left( q+1\right) <...<\sigma \left( q+r\right) }{\underset%
{\sigma \left( 1\right) <...<\sigma \left( q\right) }{\tsum }}sgn\left(
\sigma \right) \overset{q+r}{\underset{i=q+1}{\tsum }}\omega \left(
z_{\sigma \left( 1\right) },...,z_{\sigma \left( q\right) }\right) \\
\qquad \cdot \theta \left( z_{\sigma \left( q+1\right) },...,\left[
z,z_{\sigma \left( i\right) }\right] _{F,h},...,z_{\sigma \left( q+r\right)
}\right) \\
=\underset{\sigma \left( q+1\right) <...<\sigma \left( q+r\right) }{\underset%
{\sigma \left( 1\right) <...<\sigma \left( q\right) }{\tsum }}sgn\left(
\sigma \right) L_{z}\omega \left( z_{\sigma \left( 1\right) },...,\left[
z,z_{\sigma \left( i\right) }\right] _{F,h},...,z_{\sigma \left( q+r\right)
}\right) \\
\qquad \cdot \theta \left( z_{\sigma \left( q+1\right) },...,z_{\sigma
\left( q+r\right) }\right) +\underset{\sigma \left( q+1\right) <...<\sigma
\left( q+r\right) }{\underset{\sigma \left( 1\right) <...<\sigma \left(
q\right) }{\tsum }}sgn\left( \sigma \right) \overset{q+r}{\underset{i=q+1}{%
\tsum }}\omega \left( z_{\sigma \left( 1\right) },...,z_{\sigma \left(
q\right) }\right) \\
\qquad \cdot L_{z}\theta \left( z_{\sigma \left( q+1\right) },...,\left[
z,z_{\sigma \left( i\right) }\right] _{F,h},...,z_{\sigma \left( q+r\right)
}\right) \\
=\left( L_{z}\omega \wedge \theta +\omega \wedge L_{z}\theta \right) \left(
z_{1},...,z_{q+r}\right)%
\end{array}%
\end{equation*}%
it results the conclusion of the theorem. \hfill \emph{q.e.d.}

\textbf{Definition 4.6 }For any $z\in \Gamma \left( F,\nu ,N\right) $, the $%
\mathcal{F}\left( N\right) $-multilinear application
\begin{equation*}
\begin{array}{rcl}
\Lambda \left( F,\nu ,N\right) & ^{\underrightarrow{\ \ i_{z}\ \ }} &
\Lambda \left( F,\nu ,N\right) \\
\Lambda ^{q}\left( F,\nu ,N\right) \ni \omega & \longmapsto & i_{z}\omega
\in \Lambda ^{q-1}\left( F,\nu ,N\right) ,%
\end{array}%
\end{equation*}%
where
\begin{equation*}
\begin{array}{c}
i_{z}\omega \left( z_{2},...,z_{q}\right) =\omega \left(
z,z_{2},...,z_{q}\right) ,%
\end{array}%
\end{equation*}%
for any $z_{2},...,z_{q}\in \Gamma \left( F,\nu ,N\right) $, will be called
the \emph{interior product associated to the section}~$z.$\bigskip

For any $f\in \mathcal{F}\left( N\right) $, we define $\ i_{z}f=0.$

\textit{Remark 4.4}\textbf{\ }If $z\in \Gamma \left( F,\nu ,N\right) ,$ $%
\omega \in $\ $\Lambda ^{p}\left( F,\nu ,N\right) $ and $U$\ is an open
subset of $N$\ such that $z_{|U}=0$\ or $\omega _{|U}=0,$\ then $\left(
i_{z}\omega \right) _{|U}=0.$

\textbf{Theorem 4.5 }\emph{If }$z\in \Gamma \left( F,\nu ,N\right) $\emph{,
then for any }$\omega \in $\emph{\ }$\Lambda ^{q}\left( F,\nu ,N\right) $%
\emph{\ and }$\theta \in $\emph{\ }$\Lambda ^{r}\left( F,\nu ,N\right) $%
\emph{\ we obtain}
\begin{equation*}
\begin{array}{c}
i_{z}\left( \omega \wedge \theta \right) =i_{z}\omega \wedge \theta +\left(
-1\right) ^{q}\omega \wedge i_{z}\theta .%
\end{array}%
\leqno(4.7)
\end{equation*}

\emph{Proof. }Let $z_{1},...,z_{q+r}\in \Gamma \left( F,\nu ,N\right) $ be
arbitrary. We observe that
\begin{equation*}
\begin{array}{l}
i_{z_{1}}\left( \omega \wedge \theta \right) \left( z_{2},...,z_{q+r}\right)
=\left( \omega \wedge \theta \right) \left( z_{1},z_{2},...,z_{q+r}\right)
\\
=\underset{\sigma \left( q+1\right) <...<\sigma \left( q+r\right) }{\underset%
{\sigma \left( 1\right) <...<\sigma \left( q\right) }{\tsum }}sgn\left(
\sigma \right) \cdot \omega \left( z_{\sigma \left( 1\right) },...,z_{\sigma
\left( q\right) }\right) \cdot \theta \left( z_{\sigma \left( q+1\right)
},...,z_{\sigma \left( q+r\right) }\right) \\
=\underset{\sigma \left( q+1\right) <...<\sigma \left( q+r\right) }{\underset%
{1=\sigma \left( 1\right) <\sigma \left( 2\right) <...<\sigma \left(
q\right) }{\tsum }}sgn\left( \sigma \right) \cdot \omega \left(
z_{1},z_{\sigma \left( 2\right) },...,z_{\sigma \left( q\right) }\right)
\cdot \theta \left( z_{\sigma \left( q+1\right) },...,z_{\sigma \left(
q+r\right) }\right) \\
+\underset{1=\sigma \left( q+1\right) <\sigma \left( q+2\right) <...<\sigma
\left( q+r\right) }{\underset{\sigma \left( 1\right) <...<\sigma \left(
q\right) }{\tsum }}sgn\left( \sigma \right) \cdot \omega \left( z_{\sigma
\left( 1\right) },...,z_{\sigma \left( q\right) }\right) \cdot \theta \left(
z_{1},z_{\sigma \left( q+2\right) },...,z_{\sigma \left( q+r\right) }\right)
\\
=\underset{\sigma \left( q+1\right) <...<\sigma \left( q+r\right) }{\underset%
{\sigma \left( 2\right) <...<\sigma \left( q\right) }{\tsum }}sgn\left(
\sigma \right) \cdot i_{z_{1}}\omega \left( z_{\sigma \left( 2\right)
},...,z_{\sigma \left( q\right) }\right) \cdot \theta \left( z_{\sigma
\left( q+1\right) },...,z_{\sigma \left( q+r\right) }\right) \\
+\underset{\sigma \left( q+2\right) <...<\sigma \left( q+r\right) }{\underset%
{\sigma \left( 1\right) <...<\sigma \left( q\right) }{\tsum }}sgn\left(
\sigma \right) \cdot \omega \left( z_{\sigma \left( 1\right) },...,z_{\sigma
\left( q\right) }\right) \cdot i_{z_{1}}\theta \left( z_{\sigma \left(
q+2\right) },...,z_{\sigma \left( q+r\right) }\right) .%
\end{array}%
\end{equation*}

In the second sum, we have the permutation%
\begin{equation*}
\sigma =\left(
\begin{array}{ccccccc}
1 & ... & q & q+1 & q+2 & ... & q+r \\
\sigma \left( 1\right) & ... & \sigma \left( q\right) & 1 & \sigma \left(
q+2\right) & ... & \sigma \left( q+r\right)%
\end{array}%
\right) .
\end{equation*}

We observe that $\sigma =\tau \circ \tau ^{\prime }$, where%
\begin{equation*}
\tau =\left(
\begin{array}{ccccccc}
1 & 2 & ... & q+1 & q+2 & ... & q+r \\
1 & \sigma \left( 1\right) & ... & \sigma \left( q\right) & \sigma \left(
q+2\right) & ... & \sigma \left( q+r\right)%
\end{array}%
\right)
\end{equation*}%
and
\begin{equation*}
\tau ^{\prime }=\left(
\begin{array}{cccccccc}
1 & 2 & ... & q & q+1 & q+2 & ... & q+r \\
2 & 3 & ... & q+1 & 1 & q+2 & ... & q+r%
\end{array}%
\right) .
\end{equation*}

Since $\tau \left( 2\right) <...<\tau \left( q+1\right) $ and $\tau ^{\prime
}$ has $q$ inversions, it results that
\begin{equation*}
sgn\left( \sigma \right) =\left( -1\right) ^{q}\cdot sgn\left( \tau \right) .
\end{equation*}

Therefore,
\begin{equation*}
\begin{array}{l}
i_{z_{1}}\left( \omega \wedge \theta \right) \left( z_{2},...,z_{q+r}\right)
=\left( i_{z_{1}}\omega \wedge \theta \right) \left( z_{2},...,z_{q+r}\right)
\\
+\left( -1\right) ^{q}\underset{\tau \left( q+2\right) <...<\tau \left(
q+r\right) }{\underset{\tau \left( 2\right) <...<\tau \left( q\right) }{%
\tsum }}sgn\left( \tau \right) \cdot \omega \left( z_{\tau \left( 2\right)
},...,z_{\tau \left( q\right) }\right) \cdot i_{z_{1}}\theta \left( z_{\tau
\left( q+2\right) },...,z_{\tau \left( q+r\right) }\right) \\
=\left( i_{z_{1}}\omega \wedge \theta \right) \left(
z_{2},...,z_{q+r}\right) +\left( -1\right) ^{q}\left( \omega \wedge
i_{z_{1}}\theta \right) \left( z_{2},...,z_{q+r}\right) .%
\end{array}%
\end{equation*}

\hfill \emph{q.e.d.}

\textbf{Theorem 4.6 }\emph{For any }$z,v\in \Gamma \left( F,\nu ,N\right) $%
\emph{\ we obtain}%
\begin{equation*}
\begin{array}{c}
L_{v}\circ i_{z}-i_{z}\circ L_{v}=i_{\left[ z,v\right] _{F,h}}.%
\end{array}%
\leqno(4.8)
\end{equation*}

\emph{Proof.} Let $\omega \in $\emph{\ }$\Lambda ^{q}\left( F,\nu ,N\right) $
be arbitrary. Since
\begin{equation*}
\begin{array}{l}
i_{z}\left( L_{v}\omega \right) \left( z_{2},...z_{q}\right) =L_{v}\omega
\left( z,z_{2},...z_{q}\right) \\
=\Gamma \left( Th\circ \rho ,h\circ \eta \right) v\left( \omega \left(
z,z_{2},...,z_{q}\right) \right) -\omega \left( \left[ v,z\right]
_{F,h},z_{2},...,z_{q}\right) \\
-\overset{q}{\underset{i=2}{\tsum }}\omega \left( \left( z,z_{2},...,\left[
v,z_{i}\right] _{F,h},...,z_{q}\right) \right) \\
=\Gamma \left( Th\circ \rho ,h\circ \eta \right) v\left( i_{z}\omega \left(
z_{2},...,z_{q}\right) \right) -\overset{q}{\underset{i=2}{\tsum }}%
i_{z}\omega \left( z_{2},...,\left[ v,z_{i}\right] _{F,h},...,z_{q}\right)
\\
-i_{\left[ v,z\right] _{F}}\left( z_{2},...,z_{q}\right) =\left( L_{v}\left(
i_{z}\omega \right) -i_{\left[ v,z\right] _{F,h}}\right) \left(
z_{2},...,z_{q}\right) ,%
\end{array}%
\end{equation*}%
for any $z_{2},...,z_{q}\in \Gamma \left( F,\nu ,N\right) $ it result the
conclusion of the theorem.\hfill \emph{q.e.d.}

\textbf{Definition 4.7 }If $f\in \mathcal{F}\left( N\right) $ and $z\in
\Gamma \left( F,\nu ,N\right) ,$ then we define \textbf{\ }%
\begin{equation*}
\begin{array}{c}
d^{F}f\left( z\right) =\Gamma \left( Th\circ \rho ,h\circ \eta \right)
\left( z\right) f.%
\end{array}%
\end{equation*}

\textbf{Theorem 4.7 }\emph{The }$\mathcal{F}\left( N\right) $\emph{%
-multilinear application }%
\begin{equation*}
\begin{array}{c}
\begin{array}{ccc}
\Lambda ^{q}\mathbf{\ }\left( F,\nu ,N\right) & ^{\underrightarrow{\,\
d^{F}\,\ }} & \Lambda ^{q+1}\mathbf{\ }\left( F,\nu ,N\right) \\
\omega & \longmapsto & d\omega%
\end{array}%
\end{array}%
\end{equation*}%
\emph{defined by}
\begin{equation*}
\begin{array}{l}
d^{F}\omega \left( z_{0},z_{1},...,z_{q}\right) =\overset{q}{\underset{i=0}{%
\tsum }}\left( -1\right) ^{i}\Gamma \left( Th\circ \rho ,h\circ \eta \right)
z_{i}\left( \omega \left( \left( z_{0},z_{1},...,\hat{z}_{i},...,z_{q}%
\right) \right) \right) \\
~\ \ \ \ \ \ \ \ \ \ \ \ \ \ \ \ \ \ \ \ \ \ \ \ \ \ +\underset{i<j}{\tsum }%
\left( -1\right) ^{i+j}\omega \left( \left( \left[ z_{i},z_{j}\right]
_{F,h},z_{0},z_{1},...,\hat{z}_{i},...,\hat{z}_{j},...,z_{q}\right) \right) ,%
\end{array}%
\end{equation*}%
\emph{for any }$z_{0},z_{1},...,z_{q}\in \Gamma \left( F,\nu ,N\right) ,$
\emph{is unique with the following property:}%
\begin{equation*}
\begin{array}{c}
L_{z}=d^{F}\circ i_{z}+i_{z}\circ d^{F},~\forall z\in \Gamma \left( F,\nu
,N\right) .%
\end{array}%
\leqno(4.9)
\end{equation*}

This $\mathcal{F}\left( N\right) $-multilinear application will be called%
\emph{\ the exterior differentiation ope\-ra\-tor for the exterior
differential algebra of the generalized Lie algebroid }$((F,\nu
,N),[,]_{F,h},(\rho ,\eta )).$

\emph{Proof. }We verify the property $\left( 4.9\right) $ Since
\begin{equation*}
\begin{array}{l}
\left( i_{z_{0}}\circ d^{F}\right) \omega \left( z_{1},...,z_{q}\right)
=d\omega \left( z_{0},z_{1},...,z_{q}\right) \\
=\overset{q}{\underset{i=0}{\tsum }}\left( -1\right) ^{i}\Gamma \left(
Th\circ \rho ,h\circ \eta \right) z_{i}\left( \omega \left( z_{0},z_{1},...,%
\hat{z}_{i},...,z_{q}\right) \right) \\
+\underset{0\leq i<j}{\tsum }\left( -1\right) ^{i+j}\omega \left( \left[
z_{i},z_{j}\right] _{F,h},z_{0},z_{1},...,\hat{z}_{i},...,\hat{z}%
_{j},...,z_{q}\right) \\
=\Gamma \left( Th\circ \rho ,h\circ \eta \right) z_{0}\left( \omega \left(
z_{1},...,z_{q}\right) \right) \\
+\overset{q}{\underset{i=1}{\tsum }}\left( -1\right) ^{i}\Gamma \left(
Th\circ \rho ,h\circ \eta \right) z_{i}\left( \omega \left( z_{0},z_{1},...,%
\hat{z}_{i},...,z_{q}\right) \right) \\
+\underset{i=1}{\overset{q}{\tsum }}\left( -1\right) ^{i}\omega \left( \left[
z_{0},z_{i}\right] _{F,h},z_{1},...,\hat{z}_{i},...,z_{q}\right) \\
+\underset{1\leq i<j}{\tsum }\left( -1\right) ^{i+j}\omega \left( \left[
z_{i},z_{j}\right] _{F,h},z_{0},z_{1},...,\hat{z}_{i},...,\hat{z}%
_{j},...,z_{q}\right) \\
=\Gamma \left( Th\circ \rho ,h\circ \eta \right) z_{0}\left( \omega \left(
z_{1},...,z_{q}\right) \right) \\
-\underset{i=1}{\overset{q}{\tsum }}\omega \left( z_{1},...,\left[
z_{0},z_{i}\right] _{F,h},...,z_{q}\right) \\
-\overset{q}{\underset{i=1}{\tsum }}\left( -1\right) ^{i-1}\Gamma \left(
Th\circ \rho ,h\circ \eta \right) z_{i}\left( i_{z_{0}}\omega \left( \left(
z_{1},...,\hat{z}_{i},...,z_{q}\right) \right) \right) \\
-\underset{1\leq i<j}{\tsum }\left( -1\right) ^{i+j-2}i_{z_{0}}\omega \left(
\left( \left[ z_{i},z_{j}\right] _{F,h},z_{1},...,\hat{z}_{i},...,\hat{z}%
_{j},...,z_{q}\right) \right) \\
=\left( L_{z_{0}}-d^{F}\circ i_{z_{0}}\right) \omega \left(
z_{1},...,z_{q}\right) ,%
\end{array}%
\end{equation*}%
for any $z_{0},z_{1},...,z_{q}\in \Gamma \left( F,\nu ,N\right) $ it results
that the property $\left( 4.9\right) $ is satisfied.

In the following, we verify the uniqueness of the operator $d^{F}.$

Let $d^{\prime F}$ be an another exterior differentiation operator
satisfying the property $\left( 4.9\right) .$

Let $S=\left\{ q\in \mathbb{N}:d^{F}\omega =d^{\prime F}\omega ,~\forall
\omega \in \Lambda ^{q}\left( F,\nu ,N\right) \right\} $ be.

Let $z\in \Gamma \left( F,\nu ,N\right) $ be arbitrary.

We observe that $\left( 4.9\right) $ is equivalent with
\begin{equation*}
\begin{array}{c}
i_{z}\circ \left( d^{F}-d^{\prime F}\right) +\left( d^{F}-d^{\prime
F}\right) \circ i_{z}=0.%
\end{array}%
\leqno(1)
\end{equation*}

Since $i_{z}f=0,$ for any $f\in \mathcal{F}\left( N\right) ,$ it results
that
\begin{equation*}
\begin{array}{c}
\left( \left( d^{F}-d^{\prime F}\right) f\right) \left( z\right) =0,~\forall
f\in \mathcal{F}\left( N\right) .%
\end{array}%
\end{equation*}

Therefore, we obtain that
\begin{equation*}
\begin{array}{c}
0\in S.%
\end{array}%
\leqno(2)
\end{equation*}

In the following, we prove that
\begin{equation*}
\begin{array}{c}
q\in S\Longrightarrow q+1\in S%
\end{array}%
\leqno(3)
\end{equation*}

Let $\omega \in \Lambda ^{p+1}\left( F,\nu ,N\right) $ be arbitrary$.$ Since
$i_{z}\omega \in \Lambda ^{q}\left( F,\nu ,N\right) $, using the equality $%
\left( 1\right) $, it results that
\begin{equation*}
i_{z}\circ \left( d^{F}-d^{\prime F}\right) \omega =0.
\end{equation*}

We obtain that, $\left( \left( d^{F}-d^{\prime F}\right) \omega \right)
\left( z_{0},z_{1},...,z_{q}\right) =0,$ for any $z_{1},...,z_{q}\in \Gamma
\left( F,\nu ,N\right) .$

Therefore $d^{F}\omega =d^{\prime F}\omega ,$ namely $q+1\in S.$

Using the \textit{Peano's Axiom }and the affirmations $\left( 2\right) $ and
$\left( 3\right) $ it results that $S=\mathbb{N}.$

Therefore, the uniqueness is verified.\hfill \emph{q.e.d.}\medskip

Note that if $\omega =\omega _{\alpha _{1}...\alpha _{q}}t^{\alpha
_{1}}\wedge ...\wedge t^{\alpha _{q}}\in \Lambda ^{q}\left( F,\nu ,N\right) $%
, then
\begin{eqnarray*}
d^{F}\omega \left( t_{\alpha _{0}},t_{\alpha _{1}},...,t_{\alpha
_{q}}\right) &=&\overset{q}{\underset{i=0}{\tsum }}\left( -1\right)
^{i}\theta _{\alpha _{i}}^{\tilde{k}}\frac{\partial \omega _{\alpha _{0},...,%
\widehat{\alpha _{i}}...\alpha _{q}}}{\partial \varkappa ^{\tilde{k}}} \\%
[1mm]
&&+\underset{i<j}{\tsum }\left( -1\right) ^{i+j}L_{\alpha _{i}\alpha
_{j}}^{\alpha }\cdot \omega _{\alpha ,\alpha _{0},...,\widehat{\alpha _{i}}%
,...,\widehat{\alpha _{j}},...,\alpha _{q}}.
\end{eqnarray*}

Therefore, we obtain%
\begin{equation*}
\begin{array}{l}
d^{F}\omega =\left( \overset{q}{\underset{i=0}{\tsum }}\left( -1\right)
^{i}\theta _{\alpha _{i}}^{\tilde{k}}\displaystyle\frac{\partial \omega
_{\alpha _{0},...,\widehat{\alpha _{i}}...\alpha _{q}}}{\partial \varkappa ^{%
\tilde{k}}}\right. \\
\left. +\underset{i<j}{\tsum }\left( -1\right) ^{i+j}L_{\alpha _{i}\alpha
_{j}}^{\alpha }\cdot \omega _{\alpha ,\alpha _{0},...,\widehat{\alpha _{i}}%
,...,\widehat{\alpha _{j}},...,\alpha _{q}}\right) t^{\alpha _{0}}\wedge
t^{\alpha _{1}}\wedge ...\wedge t^{\alpha _{q}}.%
\end{array}%
\leqno(4.10)
\end{equation*}

\emph{Remark 4.5 }If $d^{F}$ is the exterior differentiation operator for
the generalized Lie algebroid
\begin{equation*}
\left( \left( F,\nu ,N\right) ,\left[ ,\right] _{F,h},\left( \rho ,\eta
\right) \right) ,
\end{equation*}%
$\omega \in $\emph{\ }$\Lambda ^{q}\left( F,\nu ,N\right) $\ and $U$\ is an
open subset of $N$\ such that $\omega _{|U}=0,$\ then $\left( d^{F}\omega
\right) _{|U}=0.$

\textbf{Theorem 4.8} \emph{The exterior differentiation operator }$d^{F}$%
\emph{\ given by the previous theorem has the following properties:}\medskip

\noindent 1. \emph{\ For any }$\omega \in $\emph{\ }$\Lambda ^{q}\left(
F,\nu ,N\right) $\emph{\ and }$\theta \in $\emph{\ }$\Lambda ^{r}\left(
F,\nu ,N\right) $\emph{\ we obtain }%
\begin{equation*}
\begin{array}{c}
d^{F}\left( \omega \wedge \theta \right) =d^{F}\omega \wedge \theta +\left(
-1\right) ^{q}\omega \wedge d^{F}\theta .%
\end{array}%
\leqno(4.11)
\end{equation*}

\noindent 2.\emph{\ For any }$z\in \Gamma \left( F,\nu ,N\right) $ we obtain
\begin{equation*}
\begin{array}{c}
L_{z}\circ d^{F}=d^{F}\circ L_{z}.%
\end{array}%
\leqno(4.12)
\end{equation*}

\noindent3.\emph{\ }$d^{F}\circ d^{F}=0.$

\emph{Proof.}\newline
1. Let $S=\left\{ q\in \mathbb{N}:d^{F}\left( \omega \wedge \theta \right)
=d^{F}\omega \wedge \theta +\left( -1\right) ^{q}\omega \wedge d^{F}\theta
,~\forall \omega \in \Lambda ^{q}\left( F,\nu ,N\right) \right\} $ be. Since
\begin{equation*}
\begin{array}{l}
d^{F}\left( f\wedge \theta \right) \left( z,v\right) =d^{F}\left( f\cdot
\theta \right) \left( z,v\right) \\
=\Gamma \left( Th\circ \rho ,h\circ \eta \right) z\left( f\omega \left(
v\right) \right) -\Gamma \left( Th\circ \rho ,h\circ \eta \right) v\left(
f\omega \left( z\right) \right) -f\omega \left( \left[ z,v\right]
_{F,h}\right) \\
=\Gamma \left( Th\circ \rho ,h\circ \eta \right) z\left( f\right) \cdot
\omega \left( v\right) +f\cdot \Gamma \left( Th\circ \rho ,h\circ \eta
\right) z\left( \omega \left( v\right) \right) \\
-\Gamma \left( Th\circ \rho ,h\circ \eta \right) v\left( f\right) \cdot
\omega \left( z\right) -f\cdot \Gamma \left( Th\circ \rho ,h\circ \eta
\right) v\left( \omega \left( z\right) \right) -f\omega \left( \left[ z,v%
\right] _{F,h}\right) \\
=d^{F}f\left( z\right) \cdot \omega \left( v\right) -d^{F}f\left( v\right)
\cdot \omega \left( z\right) +f\cdot d^{F}\omega \left( z,v\right) \\
=\left( d^{F}f\wedge \omega \right) \left( z,v\right) +\left( -1\right)
^{0}f\cdot d^{F}\omega \left( z,v\right) \\
=\left( d^{F}f\wedge \omega \right) \left( z,v\right) +\left( -1\right)
^{0}\left( f\wedge d^{F}\omega \right) \left( z,v\right) ,~\forall z,v\in
\Gamma \left( F,\nu ,N\right) ,%
\end{array}%
\end{equation*}%
it results that
\begin{equation*}
\begin{array}{c}
0\in S.%
\end{array}%
\leqno(1.1)
\end{equation*}

In the following we prove that
\begin{equation*}
\begin{array}{c}
q\in S\Longrightarrow q+1\in S.%
\end{array}%
\leqno(1.2)
\end{equation*}

Without restricting the generality, we consider that $\theta \in \Lambda
^{r}\left( F,\nu ,N\right) .$ Since
\begin{equation*}
\begin{array}{l}
d^{F}\left( \omega \wedge \theta \right) \left(
z_{0},z_{1},...,z_{q+r}\right) =i_{z_{0}}\circ d^{F}\left( \omega \wedge
\theta \right) \left( z_{1},...,z_{q+r}\right) \\
=L_{z_{0}}\left( \omega \wedge \theta \right) \left(
z_{1},...,z_{q+r}\right) -d^{F}\circ i_{z_{0}}\left( \omega \wedge \theta
\right) \left( z_{1},...,z_{q+r}\right) \\
=\left( L_{z_{0}}\omega \wedge \theta +\omega \wedge L_{z_{0}}\theta \right)
\left( z_{1},...,z_{q+r}\right) \\
-\left[ d^{F}\circ \left( i_{z_{0}}\omega \wedge \theta +\left( -1\right)
^{q}\omega \wedge i_{z_{0}}\theta \right) \right] \left(
z_{1},...,z_{q+r}\right) \\
=\left( L_{z_{0}}\omega \wedge \theta +\omega \wedge L_{z_{0}}\theta -\left(
d^{F}\circ i_{z_{0}}\omega \right) \wedge \theta \right) \left(
z_{1},...,z_{q+r}\right) \\
-\left( \left( -1\right) ^{q-1}i_{z_{0}}\omega \wedge d^{F}\theta +\left(
-1\right) ^{q}d^{F}\omega \wedge i_{z_{0}}\theta \right) \left(
z_{1},...,z_{q+r}\right) \\
-\left( -1\right) ^{2q}\omega \wedge d^{F}\circ i_{z_{0}}\theta \left(
z_{1},...,z_{q+r}\right) \\
=\left( \left( L_{z_{0}}\omega -d^{F}\circ i_{z_{0}}\omega \right) \wedge
\theta \right) \left( z_{1},...,z_{q+r}\right) \\
+\omega \wedge \left( L_{z_{0}}\theta -d^{F}\circ i_{z_{0}}\theta \right)
\left( z_{1},...,z_{q+r}\right) \\
+\left( \left( -1\right) ^{q}i_{z_{0}}\omega \wedge d^{F}\theta -\left(
-1\right) ^{q}d^{F}\omega \wedge i_{z_{0}}\theta \right) \left(
z_{1},...,z_{q+r}\right) \\
=\left[ \left( \left( i_{z_{0}}\circ d^{F}\right) \omega \right) \wedge
\theta +\left( -1\right) ^{q+1}d^{F}\omega \wedge i_{z_{0}}\theta \right]
\left( z_{1},...,z_{q+r}\right) \\
+\left[ \omega \wedge \left( \left( i_{z_{0}}\circ d^{F}\right) \theta
\right) +\left( -1\right) ^{q}i_{z_{0}}\omega \wedge d^{F}\theta \right]
\left( z_{1},...,z_{q+r}\right) \\
=\left[ i_{z_{0}}\left( d^{F}\omega \wedge \theta \right) +\left( -1\right)
^{q}i_{z_{0}}\left( \omega \wedge d^{F}\theta \right) \right] \left(
z_{1},...,z_{q+r}\right) \\
=\left[ d^{F}\omega \wedge \theta +\left( -1\right) ^{q}\omega \wedge
d^{F}\theta \right] \left( z_{1},...,z_{q+r}\right) ,%
\end{array}%
\end{equation*}%
for any $z_{0},z_{1},...,z_{q+r}\in \Gamma \left( F,\nu ,N\right) $, it
results $\left( 1.2\right) .$

Using the \textit{Peano's Axiom }and the affirmations $\left( 1.1\right) $
and $\left( 1.2\right) $ it results that $S=\mathbb{N}.$

Therefore, it results the conclusion of affirmation 1.\medskip

2. Let $z\in \Gamma \left( F,\nu ,N\right) $ be arbitrary.

Let $S=\left\{ q\in \mathbb{N}:\left( L_{z}\circ d^{F}\right) \omega =\left(
d^{F}\circ L_{z}\right) \omega ,~\forall \omega \in \Lambda ^{q}\left( F,\nu
,N\right) \right\} $ be.

Let $f\in \mathcal{F}\left( N\right) $ be arbitrary. Since
\begin{equation*}
\begin{array}{l}
\left( d^{F}\circ L_{z}\right) f\left( v\right) =i_{v}\circ \left(
d^{F}\circ L_{z}\right) f=\left( i_{v}\circ d^{F}\right) \circ L_{z}f \\
=\left( L_{v}\circ L_{z}\right) f-\left( \left( d^{F}\circ i_{v}\right)
\circ L_{z}\right) f \\
=\left( L_{v}\circ L_{z}\right) f-L_{\left[ z,v\right] _{F,h}}f+d^{F}\circ
i_{\left[ z,v\right] _{F,h}}f-d^{F}\circ L_{z}\left( i_{v}f\right) \\
=\left( L_{v}\circ L_{z}\right) f-L_{\left[ z,v\right] _{F,h}}f+d^{F}\circ
i_{\left[ z,v\right] _{F,h}}f-0 \\
=\left( L_{v}\circ L_{z}\right) f-L_{\left[ z,v\right] _{F,h}}f+d^{F}\circ
i_{\left[ z,v\right] _{F,h}}f-L_{z}\circ d^{F}\left( i_{v}f\right) \\
=\left( L_{z}\circ i_{v}\right) \left( d^{F}f\right) -L_{\left[ z,v\right]
_{F,h}}f+d^{F}\circ i_{\left[ z,v\right] _{F,h}}f \\
=\left( i_{v}\circ L_{z}\right) \left( d^{F}f\right) +L_{\left[ z,v\right]
_{F,h}}f-L_{\left[ z,v\right] _{F,h}}f \\
=i_{v}\circ \left( L_{z}\circ d^{F}\right) f=\left( L_{z}\circ d^{F}\right)
f\left( v\right) ,~\forall v\in \Gamma \left( F,\nu ,N\right) ,%
\end{array}%
\end{equation*}%
it results that
\begin{equation*}
\begin{array}{c}
0\in S.%
\end{array}%
\leqno(2.1)
\end{equation*}

In the following we prove that
\begin{equation*}
\begin{array}{c}
q\in S\Longrightarrow q+1\in S.%
\end{array}%
\leqno(2.2)
\end{equation*}

Let $\omega \in \Lambda ^{q}\left( F,\nu ,N\right) $ be arbitrary. Since
\begin{equation*}
\begin{array}{l}
\left( d^{F}\circ L_{z}\right) \omega \left( z_{0},z_{1},...,z_{q}\right)
=i_{z_{0}}\circ \left( d^{F}\circ L_{z}\right) \omega \left(
z_{1},...,z_{q}\right) \vspace*{1mm} \\
=\left( i_{z_{0}}\circ d^{F}\right) \circ L_{z}\omega \left(
z_{1},...,z_{q}\right) \vspace*{1mm} \\
=\left[ \left( L_{z_{0}}\circ L_{z}\right) \omega -\left( \left( d^{F}\circ
i_{z_{0}}\right) \circ L_{z}\right) \omega \right] \left(
z_{1},...,z_{q}\right) \vspace*{1mm} \\
=\left[ \left( L_{z_{0}}\circ L_{z}\right) \omega -L_{\left[ z,z_{0}\right]
_{F,h}}\omega \right] \left( z_{1},...,z_{q}\right) \vspace*{1mm} \\
+\left[ d^{F}\circ i_{\left[ z,z_{0}\right] _{F,h}}\omega -d^{F}\circ
L_{z}\left( i_{z_{0}}\omega \right) \right] \left( z_{1},...,z_{q}\right)
\vspace*{1mm} \\
\overset{ip.}{=}\left[ \left( L_{z_{0}}\circ L_{z}\right) \omega -L_{\left[
z,z_{0}\right] _{F,h}}\omega \right] \left( z_{1},...,z_{q}\right) \vspace*{%
1mm} \\
+\left[ d^{F}\circ i_{\left[ z,z_{0}\right] _{F,h}}\omega -L_{z}\circ
d^{F}\left( i_{z_{0}}\omega \right) \right] \left( z_{1},...,z_{q}\right)
\vspace*{1mm} \\
=\left[ \left( L_{z}\circ i_{z_{0}}\right) \left( d^{F}\omega \right) -L_{%
\left[ z,z_{0}\right] _{F,h}}\omega +d^{F}\circ i_{\left[ z,z_{0}\right]
_{F,h}}\omega \right] \left( z_{1},...,z_{q}\right) \vspace*{1mm} \\
=\left[ \left( i_{z_{0}}\circ L_{z}\right) \left( d^{F}\omega \right) +L_{%
\left[ z,z_{0}\right] _{F,h}}\omega -L_{\left[ z,z_{0}\right] _{F,h}}\omega %
\right] \left( z_{1},...,z_{q}\right) \vspace*{1mm} \\
=i_{z_{0}}\circ \left( L_{z}\circ d^{F}\right) \omega \left(
z_{1},...,z_{q}\right) \vspace*{1mm} \\
=\left( L_{z}\circ d^{F}\right) \omega \left( z_{0},z_{1},...,z_{q}\right)
,~\forall z_{0},z_{1},...,z_{q}\in \Gamma \left( F,\nu ,N\right) ,%
\end{array}%
\end{equation*}%
it results $\left( 2.2\right) .$

Using the \textit{Peano's Axiom }and the affirmations $\left( 2.1\right) $
and $\left( 2.2\right) $ it results that $S=\mathbb{N}.$

Therefore, it results the conclusion of affirmation 2.\medskip

3. It is remarked that
\begin{equation*}
\begin{array}{l}
i_{z}\circ \left( d^{F}\circ d^{F}\right) =\left( i_{z}\circ d^{F}\right)
\circ d^{F}=L_{z}\circ d^{F}-\left( d^{F}\circ i_{z}\right) \circ d^{F}
\vspace*{1mm} \\
=L_{z}\circ d^{F}-d^{F}\circ L_{z}+d^{F}\circ \left( d^{F}\circ i_{z}\right)
=\left( d^{F}\circ d^{F}\right) \circ i_{z},%
\end{array}%
\end{equation*}%
for any $z\in \Gamma \left( F,\nu ,N\right) .$

Let $\omega \in \Lambda ^{q}\left( F,\nu ,N\right) $ be arbitrary. Since
\begin{equation*}
\begin{array}{l}
\left( d^{F}\circ d^{F}\right) \omega \left( z_{1},...,z_{q+2}\right)
=i_{z_{q+2}}\circ ...\circ i_{z_{1}}\circ \left( d^{F}\circ d^{F}\right)
\omega =...\vspace*{1mm} \\
=i_{z_{q+2}}\circ \left( d^{F}\circ d^{F}\right) \circ i_{z_{q+1}}\left(
\omega \left( z_{1},...,z_{q}\right) \right) \vspace*{1mm} \\
=i_{z_{q+2}}\circ \left( d^{F}\circ d^{F}\right) \left( 0\right) =0,~\forall
z_{1},...,z_{q+2}\in \Gamma \left( F,\nu ,N\right) ,%
\end{array}%
\end{equation*}%
it results the conclusion of affirmation 3. \hfill \emph{q.e.d.}

\textbf{Theorem 4.9} \emph{If }$d^{F}$ \emph{is the exterior differentiation
operator for the exterior differential} $\mathcal{F}(N)$\emph{-algebra} $%
(\Lambda (F,\nu ,N),+,\cdot ,\wedge )$, \emph{then we obtain the structure
equations of Maurer-Cartan type }%
\begin{equation*}
\begin{array}{c}
d^{F}t^{\alpha }=-\displaystyle\frac{1}{2}L_{\beta \gamma }^{\alpha
}t^{\beta }\wedge t^{\gamma },~\alpha \in \overline{1,p}%
\end{array}%
\leqno(\mathcal{C}_{1})
\end{equation*}%
\emph{and\ }%
\begin{equation*}
\begin{array}{c}
d^{F}\varkappa ^{\tilde{\imath}}=\theta _{\alpha }^{\tilde{\imath}}t^{\alpha
},~\tilde{\imath}\in \overline{1,n},%
\end{array}%
\leqno(\mathcal{C}_{2})
\end{equation*}%
\emph{where }$\left\{ t^{\alpha },\alpha \in \overline{1,p}\right\} ~$\emph{%
is the coframe of the vector bundle }$\left( F,\nu ,N\right) .$\bigskip

This equations will be called \emph{the structure equations of Maurer-Cartan
type associa\-ted to the generalized Lie algebroid }$\left( \left( F,\nu
,N\right) ,\left[ ,\right] _{F,h},\left( \rho ,\eta \right) \right) .$

\emph{Proof.} Let $\alpha \in \overline{1,p}$ be arbitrary. Since%
\begin{equation*}
\begin{array}{c}
d^{F}t^{\alpha }\left( t_{\beta },t_{\gamma }\right) =-L_{\beta \gamma
}^{\alpha },~\forall \beta ,\gamma \in \overline{1,p}%
\end{array}%
\end{equation*}%
it results that
\begin{equation*}
\begin{array}{c}
d^{F}t^{\alpha }=-\underset{\beta <\gamma }{\tsum }L_{\beta \gamma }^{\alpha
}t^{\beta }\wedge t^{\gamma }.%
\end{array}%
\leqno(1)
\end{equation*}

Since $L_{\beta \gamma }^{\alpha }=-L_{\gamma \beta }^{\alpha }$ and $%
t^{\beta }\wedge t^{\gamma }=-t^{\gamma }\wedge t^{\beta }$, for nay $\beta
,\gamma \in \overline{1,p},$ it results that
\begin{equation*}
\begin{array}{c}
\underset{\beta <\gamma }{\tsum }L_{\beta \gamma }^{\alpha }t^{\beta }\wedge
t^{\gamma }=\displaystyle\frac{1}{2}L_{\beta \gamma }^{\alpha }t^{\beta
}\wedge t^{\gamma }%
\end{array}%
\leqno(2)
\end{equation*}

Using the equalities $\left( 1\right) $ and $\left( 2\right) $ it results
the structure equation $(\mathcal{C}_{1}).$

Let $\tilde{\imath}\in \overline{1,n}$ be arbitrarily. Since
\begin{equation*}
\begin{array}{c}
d^{F}\varkappa ^{\tilde{\imath}}\left( t_{\alpha }\right) =\theta _{\alpha
}^{\tilde{\imath}},~\forall \alpha \in \overline{1,p}%
\end{array}%
\end{equation*}%
it results the structure equation $(\mathcal{C}_{2}).$\hfill \emph{q.e.d.}%
\bigskip

\textbf{Corollary 4.1 }\emph{If }$\left( \left( h^{\ast }F,h^{\ast }\nu
,M\right) ,\left[ ,\right] _{h^{\ast }F},\left( \overset{h^{\ast }F}{\rho }%
,Id_{M}\right) \right) $ \emph{\ is the pull-back Lie algebroid associated
to the generalized Lie algebroid} $\left( \left( F,\nu ,N\right) ,\left[ ,%
\right] _{F,h},\left( \rho ,\eta \right) \right) $ \emph{and }$d^{h^{\ast }F}
$\emph{\ is the exterior differentiation operator for the exterior
differential }$\mathcal{F}\left( M\right) $\emph{-algebra }%
\begin{equation*}
\left( \Lambda \left( h^{\ast }F,h^{\ast }\nu ,M\right) ,+,\cdot ,\wedge
\right),
\end{equation*}%
\emph{\ then we obtain the following structure equations of Maurer-Cartan
type \ }%
\begin{equation*}
\begin{array}{c}
d^{h^{\ast }F}T^{\alpha }=-\displaystyle\frac{1}{2}\left( L_{\beta \gamma
}^{\alpha }\circ h\right) T^{\beta }\wedge T^{\gamma },~\alpha \in \overline{%
1,p}%
\end{array}%
\leqno(\mathcal{C}_{1}^{\prime })
\end{equation*}%
\emph{and\ }%
\begin{equation*}
\begin{array}{c}
d^{h^{\ast }F}x^{i}=\left( \rho _{\alpha }^{i}\circ h\right) T^{\alpha
},~i\in \overline{1,m}.%
\end{array}%
\leqno(\mathcal{C}_{2}^{\prime })
\end{equation*}

This equations will be called \emph{the structure equations of Maurer-Cartan
type associated to the pull-back Lie algebroid }%
\begin{equation*}
\left( \left( h^{\ast }F,h^{\ast }\nu ,M\right) ,\left[ ,\right] _{h^{\ast
}F},\left( \overset{h^{\ast }F}{\rho },Id_{M}\right) \right) .
\end{equation*}

\textbf{Theorem 4.10 (}of Cartan type)\textbf{\ }\emph{Let }$\left( E,\pi
,M\right) $\emph{\ be an IDS of the generalized Lie algebroid }$\left(
\left( F,\nu ,N\right) ,\left[ ,\right] _{F,h},\left( \rho ,\eta \right)
\right) .$\emph{\ If }$\left\{ \Theta ^{r+1},...,\Theta ^{p}\right\} $\emph{%
\ is a base for the }$\mathcal{F}\left( M\right) $\emph{-submodule }$\left(
\Gamma \left( E^{0},\pi ^{0},M\right) ,+,\cdot \right) $\emph{, then the IDS
}$\left( E,\pi ,M\right) $\emph{\ is involutive if and only if it exists }%
\begin{equation*}
\Omega _{\beta }^{\alpha }\in \Lambda ^{1}\left( h^{\ast }F,h^{\ast }\nu
,M\right) ,~\alpha ,\beta \in \overline{r+1,p}
\end{equation*}%
\emph{such that}
\begin{equation*}
d^{h^{\ast }F}\Theta ^{\alpha }=\Sigma _{\beta \in \overline{r+1,p}}\Omega
_{\beta }^{\alpha }\wedge \Theta ^{\beta }\in \mathcal{I}\left( \Gamma
\left( E^{0},\pi ^{0},M\right) \right) .
\end{equation*}

\emph{Proof: }Let $\left\{ S_{1},...,S_{r}\right\} $ be a base for the $%
\mathcal{F}\left( M\right) $-submodule $\left( \Gamma \left( E,\pi ,M\right)
,+,\cdot \right) $

Let $\left\{ S_{r+1},...,S_{p}\right\} \in \Gamma \left( h^{\ast }F,h^{\ast
}\nu ,M\right) $ such that
\begin{equation*}
\left\{ S_{1},...,S_{r},S_{r+1},...,S_{p}\right\}
\end{equation*}%
is a base for the $\mathcal{F}\left( M\right) $-module
\begin{equation*}
\left( \Gamma \left( h^{\ast }F,h^{\ast }\nu ,M\right) ,+,\cdot \right) .
\end{equation*}

Let $\Theta ^{1},...,\Theta ^{r}\in \Gamma \left( \overset{\ast }{h^{\ast }F}%
,\overset{\ast }{h^{\ast }\nu },M\right) $ such that
\begin{equation*}
\left\{ \Theta ^{1},...,\Theta ^{r},\Theta ^{r+1},...,\Theta ^{p}\right\}
\end{equation*}
is a base for the $\mathcal{F}\left( M\right) $-module
\begin{equation*}
\left( \Gamma \left( \overset{\ast }{h^{\ast }F},\overset{\ast }{h^{\ast
}\nu },M\right) ,+,\cdot \right) .
\end{equation*}

For any $a,b\in \overline{1,r}$ and $\alpha ,\beta \in \overline{r+1,p}$, we
have the equalities:%
\begin{equation*}
\begin{array}{ccc}
\Theta ^{a}\left( S_{b}\right) & = & \delta _{b}^{a} \\
\Theta ^{a}\left( S_{\beta }\right) & = & 0 \\
\Theta ^{\alpha }\left( S_{b}\right) & = & 0 \\
\Theta ^{\alpha }\left( S_{\beta }\right) & = & \delta _{\beta }^{\alpha }%
\end{array}%
\end{equation*}

We remark that the set of the $2$-forms%
\begin{equation*}
\left\{ \Theta ^{a}\wedge \Theta ^{b},\Theta ^{a}\wedge \Theta ^{\beta
},\Theta ^{\alpha }\wedge \Theta ^{\beta },~a,b\in \overline{1,r}\wedge
\alpha ,\beta \in \overline{r+1,p}\right\}
\end{equation*}%
is a base for the $\mathcal{F}\left( M\right) $-module
\begin{equation*}
\left( \Lambda ^{2}\left( h^{\ast }F,h^{\ast }\nu ,M\right) ,+,\cdot \right)
.
\end{equation*}

Therefore, we have%
\begin{equation*}
d^{h^{\ast }F}\Theta ^{\alpha }=\Sigma _{b<c}A_{bc}^{\alpha }\Theta
^{b}\wedge \Theta ^{c}+\Sigma _{b,\gamma }B_{b\gamma }^{\alpha }\Theta
^{b}\wedge \Theta ^{\gamma }+\Sigma _{\beta <\gamma }C_{\beta \gamma
}^{\alpha }\Theta ^{\beta }\wedge \Theta ^{\gamma },\leqno\left( 1\right)
\end{equation*}%
where, $A_{bc}^{\alpha },B_{b\gamma }^{\alpha }$ and $C_{\beta \gamma
}^{\alpha },~a,b,c\in \overline{1,r}\wedge \alpha ,\beta ,\gamma \in
\overline{r+1,p}$ are real local functions such that $A_{bc}^{\alpha
}=-A_{cb}^{\alpha }$ and $C_{\beta \gamma }^{\alpha }=-C_{\gamma \beta
}^{\alpha }.$

Using the formula%
\begin{equation*}
d^{h^{\ast }F}\Theta ^{\alpha }\left( S_{b},S_{c}\right) =\Gamma \left(
\overset{h^{\ast }F}{\rho },Id_{M}\right) S_{b}\left( \Theta ^{\alpha
}\left( S_{c}\right) \right) -\Gamma \left( \overset{h^{\ast }F}{\rho }%
,Id_{M}\right) S_{c}\left( \Theta ^{\alpha }\left( S_{b}\right) \right)
-\Theta ^{\alpha }\left( \left[ S_{b},S_{c}\right] _{h^{\ast }F}\right) ,%
\leqno\left( 2\right)
\end{equation*}%
we obtain that
\begin{equation*}
A_{bc}^{\alpha }=-\Theta ^{\alpha }\left( \left[ S_{b},S_{c}\right]
_{h^{\ast }F}\right) ,~\forall \left( b,c\in \overline{1,r}\wedge \alpha \in
\overline{r+1,p}\right) .\leqno\left( 3\right)
\end{equation*}

We admit that $\left( E,\pi ,M\right) $ is an involutive IDS of the
generalized Lie algebroid $\left( \left( F,\nu ,N\right) ,\left[ ,\right]
_{F,h},\left( \rho ,\eta \right) \right) .$

As
\begin{equation*}
\left[ S_{b},S_{c}\right] _{h^{\ast }F}\in \Gamma \left( E,\pi ,M\right)
,~\forall b,c\in \overline{1,r}
\end{equation*}%
it results that
\begin{equation*}
\Theta ^{\alpha }\left( \left[ S_{b},S_{c}\right] _{h^{\ast }F}\right)
=0,~\forall \left( b,c\in \overline{1,r}\wedge \alpha \in \overline{r+1,p}%
\right) .
\end{equation*}

Therefore,
\begin{equation*}
A_{bc}^{\alpha }=0,~\forall \left( b,c\in \overline{1,r}\wedge \alpha \in
\overline{r+1,p}\right)
\end{equation*}%
and we obtain%
\begin{equation*}
\begin{array}{ccl}
d^{h^{\ast }F}\Theta ^{\alpha } & = & \Sigma _{b,\gamma }B_{b\gamma
}^{\alpha }\Theta ^{b}\wedge \Theta ^{\gamma }+\frac{1}{2}C_{\beta \gamma
}^{\alpha }\Theta ^{\beta }\wedge \Theta ^{\gamma } \\
& = & \left( B_{b\gamma }^{\alpha }\Theta ^{b}+\frac{1}{2}C_{\beta \gamma
}^{\alpha }\Theta ^{\beta }\right) \wedge \Theta ^{\gamma }.%
\end{array}%
\end{equation*}

As
\begin{equation*}
\Omega _{\gamma }^{\alpha }\overset{put}{=}B_{b\gamma }^{\alpha }\Theta ^{b}+%
\frac{1}{2}C_{\beta \gamma }^{\alpha }\Theta ^{\beta }\in \Lambda ^{1}\left(
h^{\ast }F,h^{\ast }\nu ,M\right) ,~\forall \alpha ,\beta \in \overline{r+1,p%
}
\end{equation*}%
it results the first implication.

Conversely, we admit that it exists
\begin{equation*}
\Omega _{\beta }^{\alpha }\in \Lambda ^{1}\left( h^{\ast }F,h^{\ast }\nu
,M\right) ,~\alpha ,\beta \in \overline{r+1,p}
\end{equation*}%
such that
\begin{equation*}
d^{h^{\ast }F}\Theta ^{\alpha }=\Sigma _{\beta \in \overline{r+1,p}}\Omega
_{\beta }^{\alpha }\wedge \Theta ^{\beta },~\forall \alpha \in \overline{%
r+1,p}.\leqno\left( 4\right)
\end{equation*}

Using the affirmations $\left( 1\right) ,\left( 2\right) $ and $\left(
4\right) $ we obtain that
\begin{equation*}
A_{bc}^{\alpha }=0,~\forall \left( b,c\in \overline{1,r}\wedge \alpha \in
\overline{r+1,p}\right) .
\end{equation*}

Using the affirmation $\left( 3\right) $, we obtain
\begin{equation*}
\Theta ^{\alpha }\left( \left[ S_{b},S_{c}\right] _{h^{\ast }F}\right)
=0,~\forall \left( b,c\in \overline{1,r}\wedge \alpha \in \overline{r+1,p}%
\right) .
\end{equation*}

Therefore,
\begin{equation*}
\left[ S_{b},S_{c}\right] _{h^{\ast }F}\in \Gamma \left( E,\pi ,M\right)
,~\forall b,c\in \overline{1,r}.
\end{equation*}

Using the \emph{Proposition 3.2}, we obtain the second implication.\hfill
\emph{q.e.d.}\medskip

Let $\left( \left( F^{\prime },\nu ^{\prime },N^{\prime }\right) ,\left[ ,%
\right] _{F^{\prime },h^{\prime }},\left( \rho ^{\prime },\eta ^{\prime
}\right) \right) $ be an another generalized Lie algebroid.

\textbf{Definition 4.8} For any morphism $\left( \varphi ,\varphi
_{0}\right) $ of
\begin{equation*}
\left( \left( F,\nu ,N\right) ,\left[ ,\right] _{F,h},\left( \rho ,\eta
\right) \right)
\end{equation*}
source and
\begin{equation*}
\left( \left( F^{\prime },\nu ^{\prime },N^{\prime }\right) ,\left[ ,\right]
_{F^{\prime },h^{\prime }},\left( \rho ^{\prime },\eta ^{\prime }\right)
\right)
\end{equation*}
target we define the application
\begin{equation*}
\begin{array}{ccc}
\Lambda ^{q}\left( F^{\prime },\nu ^{\prime },N^{\prime }\right) & ^{%
\underrightarrow{\ \left( \varphi ,\varphi _{0}\right) ^{\ast }\ }} &
\Lambda ^{q}\left( F,\nu ,N\right) \\
\omega ^{\prime } & \longmapsto & \left( \varphi ,\varphi _{0}\right) ^{\ast
}\omega ^{\prime }%
\end{array}%
,
\end{equation*}%
where
\begin{equation*}
\begin{array}{c}
\left( \left( \varphi ,\varphi _{0}\right) ^{\ast }\omega ^{\prime }\right)
\left( z_{1},...,z_{q}\right) =\omega ^{\prime }\left( \Gamma \left( \varphi
,\varphi _{0}\right) \left( z_{1}\right) ,...,\Gamma \left( \varphi ,\varphi
_{0}\right) \left( z_{q}\right) \right) ,%
\end{array}%
\end{equation*}%
for any $z_{1},...,z_{q}\in \Gamma \left( F,\nu ,N\right) .$

\bigskip \textit{\noindent Remark 4.5}\textbf{\ }It is remarked that the $%
\mathbf{B}^{\mathbf{v}}$-morphism $\left( Th\circ \rho ,h\circ \eta \right) $
is a $\mathbf{GLA}$-morphism~of
\begin{equation*}
\left( \left( F,\nu ,N\right) ,\left[ ,\right] _{F,h},\left( \rho ,\eta
\right) \right)
\end{equation*}%
\emph{\ }source and
\begin{equation*}
\left( \left( TN,\tau _{N},N\right) ,\left[ ,\right] _{TN,Id_{N}},\left(
Id_{TN},Id_{N}\right) \right)
\end{equation*}%
target.

Moreover, for any $\tilde{\imath}\in \overline{1,n}$, we obtain%
\begin{equation*}
\begin{array}{c}
\left( Th\circ \rho ,h\circ \eta \right) ^{\ast }\left( d\varkappa ^{\tilde{%
\imath}}\right) =d^{F}\varkappa ^{\tilde{\imath}},%
\end{array}%
\end{equation*}%
where $d$ is the exterior differentiation operator associated to the
exterior differential Lie $\mathcal{F}\left( N\right) $-algebra\emph{\ }%
\begin{equation*}
\begin{array}{c}
\left( \Lambda \left( TN,\tau _{N},N\right) ,+,\cdot ,\wedge \right) .%
\end{array}%
\end{equation*}

\textbf{Theorem 4.11 }\emph{If }$\left( \varphi ,\varphi _{0}\right) $\emph{%
\ is a morphism of }%
\begin{equation*}
\left( \left( F,\nu ,N\right) ,\left[ ,\right] _{F,h},\left( \rho ,\eta
\right) \right)
\end{equation*}%
\emph{\ source and }%
\begin{equation*}
\left( \left( F^{\prime },\nu ^{\prime },N^{\prime }\right) ,\left[ ,\right]
_{F^{\prime },h^{\prime }},\left( \rho ^{\prime },\eta ^{\prime }\right)
\right)
\end{equation*}%
\emph{target, then the following affirmations are satisfied:}\medskip

\noindent 1. \emph{For any }$\omega ^{\prime }\in \Lambda ^{q}\left(
F^{\prime },\nu ^{\prime },N^{\prime }\right) $\emph{\ and }$\theta ^{\prime
}\in \Lambda ^{r}\left( F^{\prime },\nu ^{\prime },N^{\prime }\right) $\emph{%
\ we obtain}%
\begin{equation*}
\begin{array}{c}
\left( \varphi ,\varphi _{0}\right) ^{\ast }\left( \omega ^{\prime }\wedge
\theta ^{\prime }\right) =\left( \varphi ,\varphi _{0}\right) ^{\ast }\omega
^{\prime }\wedge \left( \varphi ,\varphi _{0}\right) ^{\ast }\theta ^{\prime
}.%
\end{array}%
\leqno(4.13)
\end{equation*}

\noindent 2.\emph{\ For any }$z\in \Gamma \left( F,\nu ,N\right) $\emph{\
and }$\omega ^{\prime }\in \Lambda ^{q}\left( F^{\prime },\nu ^{\prime
},N^{\prime }\right) $\emph{\ we obtain}%
\begin{equation*}
\begin{array}{c}
i_{z}\left( \left( \varphi ,\varphi _{0}\right) ^{\ast }\omega ^{\prime
}\right) =\left( \varphi ,\varphi _{0}\right) ^{\ast }\left( i_{\Gamma
\left( \varphi ,\varphi _{0}\right) z}\omega ^{\prime }\right) .%
\end{array}%
\leqno(4.14)
\end{equation*}

\noindent 3. \emph{If }$N=N^{\prime }$ \emph{and }%
\begin{equation*}
\left( Th\circ \rho ,h\circ \eta \right) =\left( Th^{\prime }\circ \rho
^{\prime },h^{\prime }\circ \eta ^{\prime }\right) \circ \left( \varphi
,\varphi _{0}\right) ,
\end{equation*}%
\emph{then we obtain}
\begin{equation*}
\begin{array}{c}
\left( \varphi ,\varphi _{0}\right) ^{\ast }\circ d^{F^{\prime }}=d^{F}\circ
\left( \varphi ,\varphi _{0}\right) ^{\ast }.%
\end{array}%
\leqno(4.15)
\end{equation*}

\emph{Proof.} 1. Let $\omega ^{\prime }\in \Lambda ^{q}\left( F^{\prime
},\nu ^{\prime },N^{\prime }\right) $ and $\theta ^{\prime }\in \Lambda
^{r}\left( F^{\prime },\nu ^{\prime },N^{\prime }\right) $ be arbitrary.
Since
\begin{equation*}
\begin{array}{l}
\displaystyle\left( \varphi ,\varphi _{0}\right) ^{\ast }\left( \omega
^{\prime }\wedge \theta ^{\prime }\right) \left( z_{1},...,z_{q+r}\right)
=\left( \omega ^{\prime }\wedge \theta ^{\prime }\right) \left( \Gamma
\left( \varphi ,\varphi _{0}\right) z_{1},...,\Gamma \left( \varphi ,\varphi
_{0}\right) z_{q+r}\right) \vspace*{1mm} \\
\qquad\displaystyle=\frac{1}{\left( q+r\right) !}\underset{\sigma \in \Sigma
_{q+r}}{\tsum }sgn\left( \sigma \right) \cdot \omega ^{\prime }\left( \Gamma
\left( \varphi ,\varphi _{0}\right) z_{1},...,\Gamma \left( \varphi ,\varphi
_{0}\right) z_{q}\right) \vspace*{1mm} \\
\hfill\cdot \theta ^{\prime }\left( \Gamma \left( \varphi ,\varphi
_{0}\right) z_{q+1},...,\Gamma \left( \varphi ,\varphi _{0}\right)
z_{q+r}\right) \vspace*{1mm} \\
\qquad\displaystyle=\frac{1}{\left( q+r\right) !}\underset{\sigma \in \Sigma
_{q+r}}{\tsum }sgn\left( \sigma \right) \cdot \left( \varphi ,\varphi
_{0}\right) ^{\ast }\omega ^{\prime }\left( z_{1},...,z_{q}\right) \left(
\varphi ,\varphi _{0}\right) ^{\ast }\theta ^{\prime }\left(
z_{q+1},...,z_{q+r}\right) \vspace*{1mm} \\
\qquad\displaystyle=\left( \left( \varphi ,\varphi _{0}\right) ^{\ast
}\omega ^{\prime }\wedge \left( \varphi ,\varphi _{0}\right) ^{\ast }\theta
^{\prime }\right) \left( z_{1},...,z_{q+r}\right) ,%
\end{array}%
\end{equation*}%
for any $z_{1},...,z_{q+r}\in \Gamma \left( F,\nu ,N\right) $, it results
the conclusion of affirmation 1.\medskip

2. Let $z\in \Gamma \left( F,\nu ,N\right) $\emph{\ }and\emph{\ }$\omega
^{\prime }\in \Lambda ^{q}\left( F^{\prime },\nu ^{\prime },N^{\prime
}\right) $ be arbitrary. Since
\begin{equation*}
\begin{array}{ll}
i_{z}\left( \left( \varphi ,\varphi _{0}\right) ^{\ast }\omega ^{\prime
}\right) \left( z_{2},...,z_{q}\right) & \displaystyle =\omega ^{\prime
}\left( \Gamma \left( \varphi ,\varphi _{0}\right) z,\Gamma \left( \varphi
,\varphi _{0}\right) z_{2},...,\Gamma \left( \varphi ,\varphi _{0}\right)
z_{q}\right) \vspace*{1mm} \\
& \displaystyle=i_{\Gamma \left( \varphi ,\varphi _{0}\right) z}\omega
^{\prime }\left( \Gamma \left( \varphi ,\varphi _{0}\right) z_{2},...,\Gamma
\left( \varphi ,\varphi _{0}\right) z_{q}\right) \vspace*{1mm} \\
& \displaystyle=\left( \varphi ,\varphi _{0}\right) ^{\ast }\left( i_{\Gamma
\left( \varphi ,\varphi _{0}\right) z}\omega ^{\prime }\right) \left(
z_{2},...,z_{q}\right) ,%
\end{array}%
\end{equation*}%
for any $z_{2},...,z_{q}\in \Gamma \left( F,\nu ,N\right) $, it results the
conclusion of affirmation 2$.$\smallskip

3. Let $\omega ^{\prime }\in \Lambda ^{q}\left( F^{\prime },\nu ^{\prime
},N^{\prime }\right) $ and $z_{0},...,z_{q}\in \Gamma \left( F,\nu ,N\right)
$ be arbitrary. Since
\begin{equation*}
\begin{array}{l}
\left( \left( \varphi ,\varphi _{0}\right) ^{\ast }d^{F^{\prime }}\omega
^{\prime }\right) \left( z_{0},...,z_{q}\right) =\left( d^{F^{\prime
}}\omega ^{\prime }\right) \left( \Gamma \left( \varphi ,\varphi _{0}\right)
z_{0},...,\Gamma \left( \varphi ,\varphi _{0}\right) z_{q}\right) \vspace*{%
1mm} \\
=\overset{q}{\underset{i=0}{\tsum }}\left( -1\right) ^{i}\Gamma \left(
Th^{\prime }\circ \rho ^{\prime },h^{\prime }\circ \eta ^{\prime }\right)
\left( \Gamma \left( \varphi ,\varphi _{0}\right) z_{i}\right) \vspace*{1mm}
\\
\hfill \cdot \omega ^{\prime }\left( \left( \Gamma \left( \varphi ,\varphi
_{0}\right) z_{0},\Gamma \left( \varphi ,\varphi _{0}\right) z_{1},...,%
\widehat{\Gamma \left( \varphi ,\varphi _{0}\right) z_{i}},...,\Gamma \left(
\varphi ,\varphi _{0}\right) z_{q}\right) \right) \vspace*{1mm} \\
+\underset{0\leq i<j}{\tsum }\left( -1\right) ^{i+j}\cdot \omega ^{\prime
}\left( \Gamma \left( \varphi ,\varphi _{0}\right) \left[ z_{i},z_{j}\right]
_{F,h},\Gamma \left( \varphi ,\varphi _{0}\right) z_{0},\Gamma \left(
\varphi ,\varphi _{0}\right) z_{1},...,\right. \vspace*{1mm} \\
\hfill \left. \cdot \widehat{\Gamma \left( \varphi ,\varphi _{0}\right) z_{i}%
},...,\widehat{\Gamma \left( \varphi ,\varphi _{0}\right) z_{j}},...,\Gamma
\left( \varphi ,\varphi _{0}\right) z_{q}\right)%
\end{array}%
\end{equation*}%
and
\begin{equation*}
\begin{array}{l}
d^{F}\left( \left( \varphi ,\varphi _{0}\right) ^{\ast }\omega ^{\prime
}\right) \left( z_{0},...,z_{q}\right) \vspace*{1mm} \\
=\overset{q}{\underset{i=0}{\tsum }}\left( -1\right) ^{i}\Gamma \left(
Th\circ \rho ,h\circ \eta \right) \left( z_{i}\right) \cdot \left( \left(
\varphi ,\varphi _{0}\right) ^{\ast }\omega ^{\prime }\right) \left(
z_{0},...,\widehat{z_{i}},...,z_{q}\right) \vspace*{1mm} \\
+\underset{0\leq i<j}{\tsum }\left( -1\right) ^{i+j}\cdot \left( \left(
\varphi ,\varphi _{0}\right) ^{\ast }\omega ^{\prime }\right) \left( \left[
z_{i},z_{j}\right] _{F,h},z_{0},...,\widehat{z_{i}},...,\widehat{z_{j}}%
,...,z_{q}\right) \vspace*{1mm} \\
=\overset{q}{\underset{i=0}{\tsum }}\left( -1\right) ^{i}\Gamma \left(
Th\circ \rho ,h\circ \eta \right) \left( z_{i}\right) \cdot \omega ^{\prime
}\left( \Gamma \left( \varphi ,\varphi _{0}\right) z_{0},...,\widehat{\Gamma
\left( \varphi ,\varphi _{0}\right) z_{i}},...,\Gamma \left( \varphi
,\varphi _{0}\right) z_{q}\right) \vspace*{1mm} \\
+\underset{0\leq i<j}{\tsum }\left( -1\right) ^{i+j}\cdot \omega ^{\prime
}\left( \Gamma \left( \varphi ,\varphi _{0}\right) \left[ z_{i},z_{j}\right]
_{F,h},\Gamma \left( \varphi ,\varphi _{0}\right) z_{0},\Gamma \left(
\varphi ,\varphi _{0}\right) z_{1},...,\right. \vspace*{1mm} \\
\left. \widehat{\Gamma \left( \varphi ,\varphi _{0}\right) z_{i}},...,%
\widehat{\Gamma \left( \varphi ,\varphi _{0}\right) z_{j}},...,\Gamma \left(
\varphi ,\varphi _{0}\right) z_{q}\right)%
\end{array}%
\end{equation*}%
it results the conclusion of affirmation 3. \hfill \emph{q.e.d.}

\textbf{Definition 4.9} For any $q\in \overline{1,n}$ we define%
\begin{equation*}
\begin{array}{c}
\mathcal{Z}^{q}\left( F,\nu ,N\right) =\left\{ \omega \in \Lambda ^{q}\left(
F,\nu ,N\right) :d\omega =0\right\} ,%
\end{array}%
\end{equation*}%
the set of \emph{closed differential exterior }$q$\emph{-forms} and
\begin{equation*}
\begin{array}{c}
\mathcal{B}^{q}\left( F,\nu ,N\right) =\left\{ \omega \in \Lambda ^{q}\left(
F,\nu ,N\right) :\exists \eta \in \Lambda ^{q-1}\left( F,\nu ,N\right)
~|~d\eta =\omega \right\} ,%
\end{array}%
\end{equation*}%
the set of \emph{exact differential exterior }$q$\emph{-forms}.

\section{Exterior Differential Systems}

Let $\left( \left( h^{\ast }F,h^{\ast }\nu ,M\right) ,\left[ ,\right]
_{h^{\ast }F},\left( \overset{h^{\ast }F}{\rho },Id_{M}\right) \right) $ be
the pull-back Lie algebroid of the generalized Lie algebroid $\left( \left(
F,\nu ,N\right) ,\left[ ,\right] _{F,h},\left( \rho ,\eta \right) \right) $.

\textbf{Definition 5.1 }Any ideal $\left( \mathcal{I},+,\cdot \right) $ of
the exterior differential algebra of the pull-back Lie algebroid $\left(
\left( h^{\ast }F,h^{\ast }\nu ,M\right) ,\left[ ,\right] _{h^{\ast
}F},\left( \overset{h^{\ast }F}{\rho },Id_{M}\right) \right) $ closed under
differentiation operator $d^{h^{\ast }F},$ namely $d^{h^{\ast }F}\mathcal{%
I\subseteq I},$ will be called \emph{differential ideal of the generalized
Lie algebroid }$\left( \left( F,\nu ,N\right) ,\left[ ,\right] _{F,h},\left(
\rho ,\eta \right) \right) .$

In particular, if $h=Id_{N}=\eta $, then we obtain the definition of the
differential ideal of a Lie algebroid.(see$\left[ 2\right] $)

\textbf{Definition 5.2 }Let $\left( \mathcal{I},+,\cdot \right) $ be a
differential ideal of the generalized Lie algebroid $\left( \left( F,\nu
,N\right) ,\left[ ,\right] _{F,h},\left( \rho ,\eta \right) \right) $.

If it exists an IDS $\left( E,\pi ,M\right) $ such that for all $k\in
\mathbb{N}^{\ast }$ and $\omega \in \mathcal{I}\cap \Lambda ^{k}\left(
h^{\ast }F,h^{\ast }\nu ,M\right) $ we have $\omega \left(
u_{1},...,u_{k}\right) =0,$ for any $u_{1},...,u_{k}\in \Gamma \left( E,\pi
,M\right) ,$ then we will say that $\left( \mathcal{I},+,\cdot \right) $%
\emph{\ is an exterior differential system (EDS) of the generalized Lie
algebroid }%
\begin{equation*}
\left( \left( F,\nu ,N\right) ,\left[ ,\right] _{F,h},\left( \rho ,\eta
\right) \right) .
\end{equation*}

In particular, if $h=Id_{N}=\eta $, then we obtain the defintion of the EDS
of a Lie algebroid.(see$\left[ 2\right] $)

\textbf{Theorem 5.1 (}of Cartan type) \emph{The IDS }$\left( E,\pi ,M\right)
$\emph{\ of the generalized Lie algebroid }$\left( \left( F,\nu ,N\right) ,%
\left[ ,\right] _{F,h},\left( \rho ,\eta \right) \right) $\emph{\ is
involutive, if and only if the ideal generated by the }$\mathcal{F}\left(
M\right) $\emph{-submodule }$\left( \Gamma \left( E^{0},\pi ^{0},M\right)
,+,\cdot \right) $\emph{\ is an EDS of the same generalized Lie algebroid.}

\emph{Proof. }Let $\left( E,\pi ,M\right) $ be an involutive IDS of the
generalized Lie algebroid
\begin{equation*}
\left( \left( F,\nu ,N\right) ,\left[ ,\right] _{F,h},\left( \rho ,\eta
\right) \right) .
\end{equation*}

Let $\left\{ \Theta ^{r+1},...,\Theta ^{p}\right\} $ be a base for the $%
\mathcal{F}\left( M\right) $-submodule $\left( \Gamma \left( E^{0},\pi
^{0},M\right) ,+,\cdot \right) .$

We know that
\begin{equation*}
\mathcal{I}\left( \Gamma \left( E^{0},\pi ^{0},M\right) \right) =\cup _{q\in
\mathbb{N}}\left\{ \Omega _{\alpha }\wedge \Theta ^{\alpha },~\left\{ \Omega
_{r+1},...,\Omega _{p}\right\} \subset \Lambda ^{q}\left( h^{\ast }F,h^{\ast
}\nu ,M\right) \right\} .
\end{equation*}

Let $q\in \mathbb{N}$ and $\left\{ \Omega _{r+1},...,\Omega _{p}\right\}
\subset \Lambda ^{q}\left( h^{\ast }F,h^{\ast }\nu ,M\right) $ be arbitrary.

Using the \emph{Theorems 4.8 and 4.10} we obtain
\begin{equation*}
\begin{array}{ccl}
d^{h^{\ast }F}\left( \Omega _{\alpha }\wedge \Theta ^{\alpha }\right) & = &
d^{h^{\ast }F}\Omega _{\alpha }\wedge \Theta ^{\alpha }+\left( -1\right)
^{q+1}\Omega _{\beta }\wedge d^{h^{\ast }F}\Theta ^{\beta } \\
& = & \left( d^{h^{\ast }F}\Omega _{\alpha }+\left( -1\right) ^{q+1}\Omega
_{\beta }\wedge \Omega _{\alpha }^{\beta }\right) \wedge \Theta ^{\alpha }.%
\end{array}%
\end{equation*}

As
\begin{equation*}
d^{h^{\ast }F}\Omega _{\alpha }+\left( -1\right) ^{q+1}\Omega _{\beta
}\wedge \Omega _{\alpha }^{\beta }\in \Lambda ^{q+2}\left( h^{\ast
}F,h^{\ast }\nu ,M\right)
\end{equation*}%
it results that
\begin{equation*}
d^{h^{\ast }F}\left( \Omega _{\beta }\wedge \Theta ^{\beta }\right) \in
\mathcal{I}\left( \Gamma \left( E^{0},\pi ^{0},M\right) \right)
\end{equation*}

Therefore,
\begin{equation*}
d^{h^{\ast }F}\mathcal{I}\left( \Gamma \left( E^{0},\pi ^{0},M\right)
\right) \subseteq \mathcal{I}\left( \Gamma \left( E^{0},\pi ^{0},M\right)
\right) .
\end{equation*}

Conversely, let $\left( E,\pi ,M\right) $ be an IDS of the generalized Lie
algebroid
\begin{equation*}
\left( \left( F,\nu ,N\right) ,\left[ ,\right] _{F,h},\left( \rho ,\eta
\right) \right)
\end{equation*}%
such that the $\mathcal{F}\left( M\right) $-submodule $\left( \mathcal{I}%
\left( \Gamma \left( E^{0},\pi ^{0},M\right) \right) ,+,\cdot \right) $ is
an EDS of the generalized Lie algebroid $\left( \left( F,\nu ,N\right) ,%
\left[ ,\right] _{F,h},\left( \rho ,\eta \right) \right) .$

Let $\left\{ \Theta ^{r+1},...,\Theta ^{p}\right\} $ be a base for the $%
\mathcal{F}\left( M\right) $-submodule $\left( \Gamma \left( E^{0},\pi
^{0},M\right) ,+,\cdot \right) .$ As
\begin{equation*}
d^{h^{\ast }F}\mathcal{I}\left( \Gamma \left( E^{0},\pi ^{0},M\right)
\right) \subseteq \mathcal{I}\left( \Gamma \left( E^{0},\pi ^{0},M\right)
\right)
\end{equation*}%
it results that it exists
\begin{equation*}
\Omega _{\beta }^{\alpha }\in \Lambda ^{1}\left( h^{\ast }F,h^{\ast }\nu
,M\right) ,~\alpha ,\beta \in \overline{r+1,p}
\end{equation*}%
such that
\begin{equation*}
d^{h^{\ast }F}\Theta ^{\alpha }=\Sigma _{\beta \in \overline{r+1,p}}\Omega
_{\beta }^{\alpha }\wedge \Theta ^{\beta }\in \mathcal{I}\left( \Gamma
\left( E^{0},\pi ^{0},M\right) \right) .
\end{equation*}

Using the \emph{Theorem 4.10}, it results that $\left( E,\pi ,M\right) $ is
an involutive IDS.\hfill \emph{q.e.d.}\medskip

\section{Torsion and curvature forms. Identities of Cartan and Bianchi type}

Using the theory of linear connections of Eresmann type presented in $\left[
1\right] $ for the diagram:%
\begin{equation*}
\begin{array}{c}
\xymatrix{E\ar[d]_\pi&\left( F,\left[ , \right] _{F,h},\left( \rho
,Id_{N}\right) \right)\ar[d]^\nu \\ M\ar[r]^h&N}%
\end{array}%
\leqno(6.1)
\end{equation*}%
where $\left( E,\pi ,M\right) \in \left\vert \mathbf{B}^{\mathbf{v}%
}\right\vert $ and $\left( \left( F,\nu ,N\right) ,\left[ ,\right]
_{F,h},\left( \rho ,Id_{N}\right) \right) \in \left\vert \mathbf{GLA}%
\right\vert ,$ we obtain a linear $\rho $-connection $\rho \Gamma $ for the
vector bundle $\left( E,\pi ,M\right) $ by components $\rho \Gamma _{b\alpha
}^{a}.$

Using the components of this linear $\rho $-connection, we obtain a linear $%
\rho $-connection $\rho \dot{\Gamma}$ for the vector bundle $\left( E,\pi
,M\right) $ given by the diagram:
\begin{equation*}
\begin{array}{ccl}
~\ \ \ E &  & \left( h^{\ast }F,\left[ ,\right] _{h^{\ast }F},\left( \overset%
{h^{\ast }F}{\rho },Id_{M}\right) \right) \\
\pi \downarrow &  & ~\ \ \downarrow h^{\ast }\nu \\
~\ \ \ M & ^{\underrightarrow{~\ \ \ \ Id_{M}~\ \ }} & ~\ \ \ M%
\end{array}%
\leqno(6.2)
\end{equation*}

If $\left( E,\pi ,M\right) =\left( F,\nu ,N\right) ,$ then, using the
components of the same linear $\rho $-connection $\rho \Gamma ,$ we can
consider a linear $\rho $-connection $\rho \ddot{\Gamma}$ for the vector
bundle $\left( h^{\ast }E,h^{\ast }\pi ,M\right) $ given by the diagram:
\begin{equation*}
\begin{array}{ccl}
~\ \ ~\ \ \ \ h^{\ast }E &  & \left( h^{\ast }E,\left[ ,\right] _{h^{\ast
}E},\left( \overset{h^{\ast }E}{\rho },Id_{M}\right) \right) \\
h^{\ast }\pi \downarrow &  & ~\ \ \downarrow h^{\ast }\pi \\
~\ \ \ \ \ M & ^{\underrightarrow{~\ \ \ \ Id_{M}~\ \ }} & ~\ \ M%
\end{array}%
,\leqno(6.3)
\end{equation*}%
\bigskip \noindent

\textbf{Definition 6.1 }If $\left( E,\pi ,M\right) =\left( F,\nu ,N\right) $%
, then the application
\begin{equation*}
\begin{array}{ccc}
\Gamma \left( h^{\ast }E,h^{\ast }\pi ,M\right) ^{2} & ^{\underrightarrow{\
\left( \rho ,h\right) \mathbb{T}\ }} & \Gamma \left( h^{\ast }E,h^{\ast }\pi
,M\right) \\
\left( U,V\right) & \longrightarrow & \rho \mathbb{T}\left( U,V\right)%
\end{array}%
\leqno(6.4)
\end{equation*}%
defined by:
\begin{equation*}
\begin{array}{c}
\left( \rho ,h\right) \mathbb{T}\left( U,V\right) =\rho \ddot{D}_{U}V-\rho
\ddot{D}_{V}U-\left[ U,V\right] _{h^{\ast }E},\,%
\end{array}%
\leqno(6.5)
\end{equation*}%
for any $U,V\in \Gamma \left( h^{\ast }E,h^{\ast }\pi ,M\right) ,$ will be
called $\left( \rho ,h\right) $\emph{-torsion associated to linear }$\rho $%
\emph{-connection }$\rho \Gamma .$

\bigskip \textit{\noindent Remark 6.1}\textbf{\ }In particular,\ if $%
h=Id_{M} $, then we obtain the application
\begin{equation*}
\begin{array}{ccc}
\Gamma \left( E,\pi ,M\right) ^{2} & ^{\underrightarrow{\ \rho \mathbb{T}\ }}
& \Gamma \left( E,\pi ,M\right) \\
\left( u,v\right) & \longrightarrow & \rho \mathbb{T}\left( u,v\right)%
\end{array}%
\leqno(6.4^{\prime })
\end{equation*}%
defined by:
\begin{equation*}
\begin{array}{c}
\rho \mathbb{T}\left( u,v\right) =\rho D_{u}v-\rho D_{v}u-\left[ u,v\right]
_{E},\,%
\end{array}%
\leqno(6.5^{\prime })
\end{equation*}%
for any $u,v\in \Gamma \left( E,\pi ,M\right) ,$\ which will be called $\rho
$\emph{-torsion associated to linear }$\rho $\emph{-connection }$\rho \Gamma
.$

Moreover, if $\rho =Id_{TM}$, then we obtain the torsion $\mathbb{T}$
associated to linear connection~$\Gamma .$

\textbf{Proposition 6.1 }\emph{The }$\left( \rho ,h\right) $\emph{-torsion }$%
\left( \rho ,h\right) \mathbb{T}$\emph{\ associated to linear }$\rho $\emph{%
-connection }$\rho \Gamma $\emph{\ is }$\mathbb{R}$\emph{-bilinear and
antisymmetric.}

\emph{If }%
\begin{equation*}
\left( \rho ,h\right) \mathbb{T}\left( S_{a},S_{b}\right) \overset{put}{=}%
\left( \rho ,h\right) \mathbb{T}_{~ab}^{c}S_{c}
\end{equation*}%
\emph{\ then }%
\begin{equation*}
\begin{array}{c}
\left( \rho ,h\right) \mathbb{T}_{~ab}^{c}=\rho \Gamma _{ab}^{c}-\rho \Gamma
_{ba}^{c}-L_{ab}^{c}\circ h.%
\end{array}%
\leqno(6.6)
\end{equation*}

\emph{In particular, if }$h=Id_{M}$\emph{\ and }$\rho \mathbb{T}\left(
s_{a},s_{b}\right) \overset{put}{=}\rho \mathbb{T}_{ab}^{c}s_{c}$\emph{,
then }%
\begin{equation*}
\begin{array}{c}
\rho \mathbb{T}_{~ab}^{c}=\rho \Gamma _{ab}^{c}-\rho \Gamma
_{ba}^{c}-L_{ab}^{c}.%
\end{array}%
\leqno(6.6^{\prime })
\end{equation*}

\emph{Moreover, if }$\rho =Id_{TM}$\emph{, then the equality }$\left(
6.6^{\prime }\right) $\emph{\ becomes:}
\begin{equation*}
\begin{array}{c}
\mathbb{T}_{~jk}^{i}=\Gamma _{jk}^{i}-\Gamma _{kj}^{i}.%
\end{array}%
\leqno(6.6^{\prime \prime })
\end{equation*}

\textbf{Definition 6.2 }If $\left( E,\pi ,M\right) =\left( F,\nu ,N\right) $%
, then the vector valued 2-form%
\begin{equation*}
\begin{array}{c}
\left( \rho ,h\right) \mathbb{T}=\left( \left( \rho ,h\right) \mathbb{T}%
_{~ab}^{c}S_{c}\right) S^{a}\wedge S^{b}%
\end{array}%
\leqno(6.7)
\end{equation*}%
will be called the \emph{vector valued form of }$\left( \rho ,h\right) $%
\emph{-torsion }$\left( \rho ,h\right) \mathbb{T}$\emph{.}

In particular, if $h=Id_{M}$, then the vector valued 2-form
\begin{equation*}
\begin{array}{c}
\rho \mathbb{T}=\left( \rho \mathbb{T}_{~ab}^{c}s_{c}\right) s^{a}\wedge
s^{b}%
\end{array}%
\leqno(6.7^{\prime })
\end{equation*}%
will be called the \emph{vector form of }$\rho $\emph{-torsion }$\rho
\mathbb{T}$\emph{.}

Moreover, if $\rho =Id_{TM}$, then the vector valued form $\left(
6.7^{\prime }\right) $ becomes:
\begin{equation*}
\begin{array}{c}
\mathbb{T}=\left( \mathbb{T}_{~jk}^{i}\frac{\partial }{\partial x^{i}}%
\right) dx^{j}\wedge dx^{k}.%
\end{array}%
\leqno(6.7^{\prime \prime })
\end{equation*}

\textbf{Definition 6.3 }For each $c\in \overline{1,n}$ we obtain the \emph{%
scalar }$2$\emph{-form of }$\left( \rho ,h\right) $\emph{-torsion }$\left(
\rho ,h\right) \mathbb{T}$
\begin{equation*}
\begin{array}{c}
\left( \rho ,h\right) \mathbb{T}^{c}=\left( \rho ,h\right) \mathbb{T}%
_{~ab}^{c}S^{a}\wedge S^{b}.%
\end{array}%
\leqno(6.8)
\end{equation*}

In particular, if $h=Id_{M}$, then, for each $c\in \overline{1,n},$ we
obtain the \emph{scalar }$2$\emph{-form of }$\rho $\emph{-torsion }$\rho
\mathbb{T}$
\begin{equation*}
\begin{array}{c}
\rho \mathbb{T}^{c}=\rho \mathbb{T}_{~ab}^{c}s^{a}\wedge s^{b}.%
\end{array}%
\leqno(6.8^{\prime })
\end{equation*}

Moreover, if $\rho =Id_{TM}$, then the scalar $2$-form $\left( 6.9^{\prime
}\right) $ becomes:
\begin{equation*}
\begin{array}{c}
\mathbb{T}^{i}=\mathbb{T}_{~jk}^{i}dx^{j}\wedge dx^{k}.%
\end{array}%
\leqno(6.8^{\prime \prime })
\end{equation*}

\textbf{Definition 6.4 }The application
\begin{equation*}
\begin{array}{ccl}
(\Gamma \left( h^{\ast }F,h^{\ast }\nu ,M\right) ^{2}{\times }\Gamma (E,\pi
,M) & ^{\underrightarrow{\ \left( \rho ,h\right) \mathbb{R}\ }} & \Gamma
(E,\pi ,M) \\
((Z,V),u) & \longrightarrow & \rho \mathbb{R}(Z,V)u%
\end{array}%
\leqno(6.9)
\end{equation*}%
defined by
\begin{equation*}
\left( \rho ,h\right) \mathbb{R}\left( Z,V\right) u=\rho \dot{D}_{Z}\left(
\rho \dot{D}_{V}u\right) -\rho \dot{D}_{V}\left( \rho \dot{D}_{Z}u\right)
-\rho \dot{D}_{\left[ Z,V\right] _{h^{\ast }F}}u,\,\leqno(6.10)
\end{equation*}%
for any $Z,V\in \Gamma \left( h^{\ast }F,h^{\ast }\nu ,M\right) ,~u\in
\Gamma \left( E,\pi ,M\right) ,$ will be called $\left( \rho ,h\right) $%
\emph{-curvature associated to linear }$\rho $\emph{-connection }$\rho
\Gamma .$

\textit{Remark 6.2 }In particular, if $h=Id_{M}$, then we obtain the
application
\begin{equation*}
\begin{array}{ccl}
\Gamma \left( F,\nu ,M\right) ^{2}{\times }\Gamma (E,\pi ,M) & ^{%
\underrightarrow{\ \rho \mathbb{R}\ }} & \Gamma (E,\pi ,M) \\
((z,v),u) & \longrightarrow & \rho \mathbb{R}(z,v)u%
\end{array}%
\leqno(6.9^{\prime })
\end{equation*}%
defined by
\begin{equation*}
\rho \mathbb{R}\left( z,v\right) u=\rho D_{z}\left( \rho D_{v}u\right) -\rho
D_{v}\left( \rho D_{z}u\right) -\rho D_{\left[ z,v\right] _{F}}u,\,\leqno%
(6.10^{\prime })
\end{equation*}%
for any $z,v\in \Gamma \left( F,\nu ,M\right) ,~u\in \Gamma \left( E,\pi
,M\right) ,$ which will be called $\rho $\emph{-curvature associated to
linear }$\rho $\emph{-connection }$\rho \Gamma .$

Moreover, if $\rho =Id_{TM}$, then we obtain the curvature $\mathbb{R}$
associated to linear connection $\Gamma .$

\textbf{Proposition 6.2 }\emph{The }$\left( \rho ,h\right) $\emph{-curvature
}$\left( \rho ,h\right) \mathbb{R}$\emph{\ associated to linear }$\rho $%
\emph{-connection }$\rho \Gamma $\emph{, is }$\mathbb{R}$\emph{-linear in
each argument and antisymmetric in the first two arguments.}

\emph{If }%
\begin{equation*}
\left( \rho ,h\right) \mathbb{R}\left( T_{\beta },T_{\alpha }\right) s_{b}%
\overset{put}{=}\left( \rho ,h\right) \mathbb{R}_{b~\alpha \beta }^{a}s_{a},
\end{equation*}%
\emph{then }%
\begin{equation*}
\begin{array}[b]{cl}
\left( \rho ,h\right) \mathbb{R}_{b~\alpha \beta }^{a} & =\rho _{\beta
}^{j}\circ h\frac{\partial \rho \Gamma _{b\alpha }^{a}}{\partial x^{j}}+\rho
\Gamma _{e\beta }^{a}\rho \Gamma _{b\alpha }^{e}-\rho _{\alpha }^{i}\circ h%
\frac{\partial \rho \Gamma _{b\beta }^{a}}{\partial x^{i}} \\
& -\rho \Gamma _{e\alpha }^{a}\rho \Gamma _{b\beta }^{e}+\rho \Gamma
_{b\gamma }^{a}L_{\alpha \beta }^{\gamma }\circ h.%
\end{array}%
\leqno(6.11)
\end{equation*}

\emph{In particular, if }$h=Id_{M}$\emph{\ and }$\rho \mathbb{R}\left(
t_{\beta },t_{\alpha }\right) s_{b}\overset{put}{=}\rho \mathbb{R}_{b~\alpha
\beta }^{a}s_{a}$\emph{, then }%
\begin{equation*}
\begin{array}[b]{c}
\rho \mathbb{R}_{b~\alpha \beta }^{a}=\rho _{\beta }^{j}\frac{\partial \rho
\Gamma _{b\alpha }^{a}}{\partial x^{j}}+\rho \Gamma _{e\beta }^{a}\rho
\Gamma _{b\alpha }^{e}-\rho _{\alpha }^{i}\frac{\partial \rho \Gamma
_{b\beta }^{a}}{\partial x^{i}}-\rho \Gamma _{e\alpha }^{a}\rho \Gamma
_{b\beta }^{e}+\rho \Gamma _{b\gamma }^{a}L_{\alpha \beta }^{\gamma }.%
\end{array}%
\leqno(6.11^{\prime })
\end{equation*}

\emph{Moreover, if }$\rho =Id_{TM}$\emph{, then equality }$\left(
6.11^{\prime }\right) $\emph{\ becomes:}
\begin{equation*}
\begin{array}[b]{c}
\mathbb{R}_{b~hk}^{a}=\frac{\partial \Gamma _{bh}^{a}}{\partial x^{k}}%
+\Gamma _{ek}^{a}\Gamma _{bh}^{e}-\frac{\partial \Gamma _{bk}^{a}}{\partial
x^{h}}-\Gamma _{eh}^{a}\Gamma _{bk}^{e}.%
\end{array}%
\leqno(6.11^{\prime \prime })
\end{equation*}

\textbf{Definition 6.5 }The vector mixed form
\begin{equation*}
\begin{array}{c}
\left( \rho ,h\right) \mathbb{R=}\left( \left( \left( \rho ,h\right) \mathbb{%
R}_{b~\alpha \beta }^{a}s_{a}\right) T^{\alpha }\wedge T^{\beta }\right)
s^{b}%
\end{array}%
\leqno(6.12)
\end{equation*}%
will be called the \emph{vector valued form of }$\left( \rho ,h\right) $%
\emph{-curvature }$\left( \rho ,h\right) \mathbb{R}$\emph{.}

In particular, if $h=Id_{M}$, then the vector mixed form
\begin{equation*}
\begin{array}{c}
\rho \mathbb{R=}\left( \left( \rho \mathbb{R}_{b~\alpha \beta
}^{a}s_{a}\right) t^{\alpha }\wedge t^{\beta }\right) s^{b}%
\end{array}%
\leqno(6.12^{\prime })
\end{equation*}%
will be called the \emph{vector valued form of }$\rho $\emph{-curvature }$%
\rho \mathbb{R}$\emph{.}

Moreover, if $\rho =Id_{TM}$, then the vector form $\left( 6.12^{\prime
}\right) $ becomes:%
\begin{equation*}
\begin{array}{c}
\mathbb{R=}\left( \left( \mathbb{R}_{b~hk}^{a}s_{a}\right) dx^{h}\wedge
dx^{k}\right) s^{b}.%
\end{array}%
\leqno(6,12^{\prime \prime })
\end{equation*}

\textbf{Definition 6.6 }For each $a,b\in \overline{1,n}$ we obtain the \emph{%
scalar }$2$\emph{-form of }$\left( \rho ,h\right) $\emph{-curvature }$\left(
\rho ,h\right) \mathbb{R}$
\begin{equation*}
\begin{array}{c}
\left( \rho ,h\right) \mathbb{R}_{b}^{a}=\left( \rho ,h\right) \mathbb{R}%
_{b~\alpha \beta }^{a}T^{\alpha }\wedge T^{\beta }.%
\end{array}%
\leqno(6.13)
\end{equation*}

In particular, if $h=Id_{M}$, then, for each $a,b\in \overline{1,n},$ we
obtain the \emph{scalar }$2$\emph{-form of }$\rho $\emph{-curvature }$\rho
\mathbb{R}$
\begin{equation*}
\begin{array}{c}
\rho \mathbb{R}_{b}^{a}=\rho \mathbb{R}_{b~\alpha \beta }^{a}t^{\alpha
}\wedge t^{\beta }.%
\end{array}%
\leqno(6.13^{\prime })
\end{equation*}

Moreover, if $\rho =Id_{TM}$, then the scalar form $\left( 6.13^{\prime
}\right) $ becomes:
\begin{equation*}
\begin{array}{c}
\mathbb{R}_{b}^{a}=\mathbb{R}_{b~hk}^{a}dx^{h}\wedge dx^{k}.%
\end{array}%
\leqno(6.13^{\prime \prime })
\end{equation*}

\textbf{Theorem 6.1 }\emph{The identities }%
\begin{equation*}
\begin{array}{c}
\left( \rho ,h\right) \mathbb{T}^{a}=d^{h^{\ast }F}S^{a}+\Omega
_{b}^{a}\wedge S^{b},%
\end{array}%
\leqno(C_{1})
\end{equation*}%
\emph{and }%
\begin{equation*}
\begin{array}{c}
\left( \rho ,h\right) \mathbb{R}_{b}^{a}=d^{h^{\ast }F}\Omega
_{b}^{a}+\Omega _{c}^{a}\wedge \Omega _{b}^{c}%
\end{array}%
\leqno(C_{2})
\end{equation*}%
\emph{hold good. These will be called the first respectively the second
identity of Cartan type.}

\emph{Proof.} To prove the first identity we consider that $\left( E,\pi
,M\right) =\left( F,\nu ,M\right) .$ Therefore, $\Omega _{b}^{a}=\rho \Gamma
_{bc}^{a}S^{c}.$ Since
\begin{equation*}
\begin{array}{l}
d^{h^{\ast }F}S^{a}(U,V)S_{a}=((\Gamma (\overset{h^{\ast }F}{\rho }%
,Id_{M})U)S^{a}(V) \vspace*{1mm} \\
\qquad-(\Gamma (\overset{h^{\ast }F}{\rho }%
,Id_{M})V)(S^{a}(U))-S^{a}([U,V]_{h^{\ast }F}))S_{a}\vspace*{1mm} \\
\qquad=(\Gamma (\overset{h^{\ast }F}{\rho },Id_{M})U)(V^{a})-(\Gamma (%
\overset{h^{\ast }F}{\rho },Id_{M})V)(U^{a})-S^{a}([U,V]_{h^{\ast }F})S_{a}%
\vspace*{1mm} \\
\qquad=\rho \ddot{D}_{U}V-V^{b}\rho \ddot{D}_{U}S_{b}-\rho \ddot{D}%
_{V}U-U^{b}\rho \ddot{D}_{V}S_{b}-[U,V]_{h^{\ast }F}\vspace*{1mm} \\
\qquad=\left( \rho ,h\right) \mathbb{T}(U,V)-(\rho \Gamma
_{bc}^{a}V^{b}U^{c}-\rho \Gamma _{bc}^{a}U^{b}V^{c})S_{a} \vspace*{1mm} \\
\qquad=(\left( \rho ,h\right) \mathbb{T}^{a}(U,V)-\Omega _{b}^{a}\wedge
S^{b}(U,V))S_{a},%
\end{array}%
\end{equation*}%
it results the first identity.

To prove the second identity, we consider that $\left( E,\pi ,M\right) \neq
\left( F,\nu ,M\right) .$ Since
\begin{equation*}
\begin{array}{l}
\left( \rho ,h\right) \mathbb{R}_{b}^{a}\left( Z,W\right) s_{a}=\left( \rho
,h\right) \mathbb{R}\left( \left( W,Z\right) ,s_{b}\right) \vspace*{1mm} \\
\qquad=\rho \dot{D}_{Z}\left( \rho \dot{D}_{W}s_{b}\right) -\rho \dot{D}%
_{W}\left( \rho \dot{D}_{Z}s_{b}\right) -\rho \dot{D}_{\left[ Z,W\right]
_{h^{\ast }F}}s_{b} \vspace*{1mm} \\
\qquad=\rho \dot{D}_{Z}\left( \Omega _{b}^{a}\left( W\right) s_{a}\right)
-\rho \dot{D}_{W}\left( \Omega _{b}^{a}\left( Z\right) s_{a}\right) -\Omega
_{b}^{a}\left( \left[ Z,W\right] _{h^{\ast }F}\right) s_{a}\vspace*{1mm} \\
\qquad+\left( \Omega _{c}^{a}\left( Z\right) \Omega _{b}^{c}\left( W\right)
-\Omega _{c}^{a}\left( W\right) \Omega _{b}^{c}\left( Z\right) \right) s_{a}
\vspace*{1mm} \\
\qquad=\left( d^{h^{\ast }F}\Omega _{b}^{a}\left( Z,W\right) +\Omega
_{c}^{a}\wedge \Omega _{b}^{c}\left( Z,W\right) \right) s_{a}%
\end{array}%
\end{equation*}%
it results the second identity.

\textbf{Corollary 6.1 }\emph{In particular, if }$h=Id_{M}, $\emph{\ then the
identities }$(C_{1})$\emph{\ and }$(C_{2})$ \emph{become}%
\begin{equation*}
\begin{array}{c}
\rho \mathbb{T}^{a}=d^{F}s^{a}+\omega _{b}^{a}\wedge s^{b},%
\end{array}%
\leqno(C_{1}^{\prime })
\end{equation*}%
\emph{and }%
\begin{equation*}
\begin{array}{c}
\rho \mathbb{R}_{b}^{a}=d^{F}\omega _{b}^{a}+\omega _{c}^{a}\wedge \omega
_{b}^{c}%
\end{array}%
\leqno(C_{2}^{\prime })
\end{equation*}%
\emph{respectively.}

\emph{Moreover, if }$\rho =Id_{TM}$\emph{, then the identities }$\left(
C_{1}^{\prime }\right) $\emph{\ and }$\left( C_{2}^{\prime }\right) \,$\emph{%
become:}%
\begin{equation*}
\begin{array}{c}
\mathbb{T}^{i}=ddx^{i}+\omega _{j}^{i}\wedge dx^{j}=\omega _{j}^{i}\wedge
dx^{j}%
\end{array}%
\leqno(C_{1}^{\prime \prime })
\end{equation*}%
\emph{and }%
\begin{equation*}
\begin{array}{c}
\mathbb{R}_{j}^{i}=d\omega _{j}^{i}+\omega _{h}^{i}\wedge \omega _{j}^{h},%
\end{array}%
\leqno(C_{2}^{\prime \prime })
\end{equation*}%
\emph{\ respectively. }\hfill \emph{q.e.d.}

\textbf{Theorem 6.2 }\emph{The identities }%
\begin{equation*}
\begin{array}{c}
d^{h^{\ast }F}\left( \rho ,h\right) \mathbb{T}^{a}=\left( \rho ,h\right)
\mathbb{R}_{b}^{a}\wedge S^{b}-\Omega _{c}^{a}\wedge \left( \rho ,h\right)
\mathbb{T}^{c}%
\end{array}%
\leqno(B_{1})
\end{equation*}%
\emph{and }%
\begin{equation*}
\begin{array}{c}
d^{h^{\ast }F}\left( \rho ,h\right) \mathbb{R}_{b}^{a}=\left( \rho ,h\right)
\mathbb{R}_{c}^{a}\wedge \Omega _{b}^{c}-\Omega _{c}^{a}\wedge \left( \rho
,h\right) \mathbb{R}_{b}^{c},%
\end{array}%
\leqno(B_{2})
\end{equation*}%
\emph{hold good.} \emph{We will called these the first respectively the
second identity of Bianchi type.}

\emph{If the }$\left( \rho ,h\right) $\emph{-torsion is null, then the first
identity of Bianchi type becomes: }%
\begin{equation*}
\begin{array}{c}
\left( \rho ,h\right) \mathbb{R}_{b}^{a}\wedge s^{b}=0.%
\end{array}%
\leqno(\tilde{B}_{1})
\end{equation*}

\textit{Proof.} We consider $\left( E,\pi ,M\right) =\left( F,\nu ,M\right)
. $ Using the first identity of Cartan type and the equality $d^{h^{\ast
}F}\circ d^{h^{\ast }F}=0,$ we obtain:
\begin{equation*}
\begin{array}{c}
d^{h^{\ast }F}\left( \rho ,h\right) \mathbb{T}^{a}=d^{h^{\ast }F}\Omega
_{b}^{a}\wedge S^{b}-\Omega _{c}^{a}\wedge d^{h^{\ast }F}S^{c}.%
\end{array}%
\end{equation*}

Using the second identity of Cartan type and the previous identity, we
obtain:
\begin{equation*}
d^{h^{\ast }F}\left( \rho ,h\right) \mathbb{T}^{a}=\left( \left( \rho
,h\right) \mathbb{R}_{b}^{a}-\Omega _{c}^{a}\wedge \Omega _{b}^{c}\right)
\wedge S^{b}-\Omega _{c}^{a}\wedge \left( \left( \rho ,h\right) \mathbb{T}%
^{c}-\Omega _{b}^{c}\wedge S^{b}\right) .
\end{equation*}

After some calculations, we obtain the first identity of Bianchi type.

Using the second identity of Cartan type and the equality $d^{h^{\ast
}F}\circ d^{h^{\ast }F}=0,$ we obtain:
\begin{equation*}
\begin{array}{c}
d^{h^{\ast }F}\Omega _{c}^{a}\wedge \Omega _{b}^{c}-\Omega _{c}^{a}\wedge
d^{h^{\ast }F}\Omega _{b}^{c}=d^{h^{\ast }F}\left( \rho ,h\right) \mathbb{R}%
_{b}^{a}.%
\end{array}%
\end{equation*}

Using the second of Cartan type and the previous identity, we obtain:%
\begin{equation*}
\begin{array}{cc}
d^{h^{\ast }F}\left( \rho ,h\right) \mathbb{R}_{b}^{a}=\left( \left( \rho
,h\right) \mathbb{R}_{c}^{a}-\Omega _{e}^{a}\wedge \Omega _{c}^{e}\right)
\wedge \Omega _{b}^{c}-\Omega _{c}^{a}\wedge \left( \left( \rho ,h\right)
\mathbb{R}_{b}^{c}-\Omega _{e}^{c}\wedge \Omega _{b}^{e}\right) . &
\end{array}%
\end{equation*}

After some calculations, we obtain the second identity of Bianchi
type.\hfill \hfill \emph{q.e.d.}

\textbf{Corollary 6.2 }\emph{In particular, if }$h=Id_{M}, $\emph{\ then the
identities }$(B_{1})$\emph{\ and }$(B_{2})$ \emph{become }%
\begin{equation*}
\begin{array}{c}
d^{F}\rho \mathbb{T}^{a}=\rho \mathbb{R}_{b}^{a}\wedge s^{b}-\omega
_{c}^{a}\wedge \rho \mathbb{T}^{c}%
\end{array}%
\leqno(B_{1}^{\prime })
\end{equation*}%
\emph{and }%
\begin{equation*}
\begin{array}{c}
d^{F}\rho \mathbb{R}_{b}^{a}=\rho \mathbb{R}_{c}^{a}\wedge \omega
_{b}^{c}-\omega _{c}^{a}\wedge \rho \mathbb{R}_{b}^{c},%
\end{array}%
\leqno(B_{2}^{\prime })
\end{equation*}%
\emph{\ respectively.}

\emph{Moreover, if }$\rho =Id_{TM}$\emph{, then the identities }$\left(
B_{1}^{\prime }\right) $\emph{\ and }$\left( B_{2}^{\prime }\right) \,$\emph{%
become:}
\begin{equation*}
\begin{array}{c}
d\mathbb{T}^{i}=\mathbb{R}_{j}^{i}\wedge dx^{j}-\omega _{k}^{i}\wedge
\mathbb{T}^{k}%
\end{array}%
\leqno(B_{1}^{\prime \prime })
\end{equation*}%
\emph{and }%
\begin{equation*}
\begin{array}{c}
d\mathbb{R}_{j}^{i}=\mathbb{R}_{h}^{i}\wedge \omega _{j}^{h}-\omega
_{h}^{i}\wedge \mathbb{R}_{j}^{h},%
\end{array}%
\leqno(B_{2}^{\prime \prime })
\end{equation*}%
\emph{\ respectively.}





\bigskip

\bigskip

\end{document}